\documentclass[applsci,article,accept,pdftex,moreauthors]{Definitions/mdpi} 
\usepackage{amsmath}
\usepackage{graphics}
\usepackage{booktabs}

\usepackage[utf8]{inputenc}
\usepackage{textcomp}
\usepackage{chngcntr}
\usepackage{graphicx}
\usepackage{subcaption}
\usepackage{graphicx}

\usepackage{amsmath}
\usepackage[version=4]{mhchem}
\usepackage{siunitx}
\usepackage{longtable,tabularx}

\usepackage{amssymb}

\usepackage{algorithm}

\usepackage{geometry}
\usepackage{atbegshi}
\usepackage{etoolbox}
\usepackage{xcolor}
\usepackage{soul}

\usepackage{booktabs} 
\usepackage{array} 
\usepackage{caption}
\usepackage{longtable}

\usepackage{ifthen}

\usepackage{nomencl}
\usepackage{xpatch}

\usepackage{amssymb}

\newlength{\nomitemorigsep}
\firstpage{1} 
\pubvolume{1}
\issuenum{1}
\articlenumber{0}
\pubyear{2022}
\copyrightyear{2022}
\datereceived{16 November 2021} 
\dateaccepted{15 December 2021} 
\datepublished{26 December 2021} 
\hreflink{https://doi.org/10.3390/app12010207} 

\Title{Benefits of Advanced Air Mobility for Society and Environment: A Case Study of Ohio}

\TitleCitation{Benefits of Advanced Air Mobility for Society and Environment: A Case Study of Ohio}

\Author{Esrat F. Dulia 
 $^{1,}$*, Mir S. Sabuj $^{2}$, and Syed A. M. Shihab 
$^{1}$}

\AuthorNames{Esrat F. Dulia, Mir S. Sabuj and Syed A. M. Shihab}

\AuthorCitation{Dulia, E.F.; Sabuj, M.S.; Shihab, S.A.M.}

\address{%
$^{1}$ \quad College of Aeronautics and Engineering, Kent State University, Kent, OH 44242, USA; sshihab@kent.edu\\
$^{2}$ \quad Pathao Ltd., Dhaka 1213, Bangladesh; shahriar052@gmail.com\\}

\corres{Correspondence: edulia@kent.edu}

\abstract{Advanced Air Mobility (AAM) is an emerging transportation system that will enable the safe and efficient low altitude operations and applications of unmanned aircraft (e.g., passenger transportation and cargo delivery) in the national airspace. This system is currently under active research and development by NASA in collaboration with FAA, other federal partner agencies, industry, and academia to develop its infrastructure, information architecture, software functions, concepts of operation, operations management tools and other functional components. Existing studies have, however, not thoroughly analyzed the net positive impact of AAM on society and environment to justify investments in its infrastructure and implementation. In this work, we fill this gap by evaluating the non-monetary social impact of AAM in the state of Ohio for passengers, patients, farmers, logistics companies and their customers and bridge inspection entities, as well as its environmental impact, by conducting a thorough data-driven quantitative cost–benefit analysis of AAM from the perspective of the state government. To this end, the most relevant and significant benefit and cost factors are identified, monetized, and estimated. Existing ground transportation for the movement of passengers and goods within and across urban areas is considered as the base case. The findings demonstrate that AAM's benefits are large and varied, far outweighing its costs. Insights on these benefits can help gain community acceptance of AAM, which is critical for successful implementation of AAM. The findings support decision-making for policymakers and provide justification for investments in AAM infrastructure by the government and private sector.}


\begin{document}

\section{Introduction} \label{intro}

AAM infrastructure will enable and support a wide range of applications, such as passenger transportation, cargo delivery, infrastructure inspection, and precision agriculture by uncrewed aerial systems (UAS) in the low altitude National Airspace System (NAS) \cite{Prevot2016}. The number of daily operations needed to fulfill the demands of these use cases is projected to approach millions in the coming years, outstripping the capacity of the traditional human-centered NAS \cite{FAA2020}. Therefore, a digital and flexible transportation framework is being developed that can coordinate among the new NAS users (the AAM users) and existing NAS users and manage the uncrewed air traffic to support ongoing efforts to integrate uncrewed aircraft into the lower airspace and ensure their safety and security from environmental risks. This framework is referred to as the AAM architecture. 

\begin{figure}[H]
    \includegraphics[width=13cm,height=8cm]{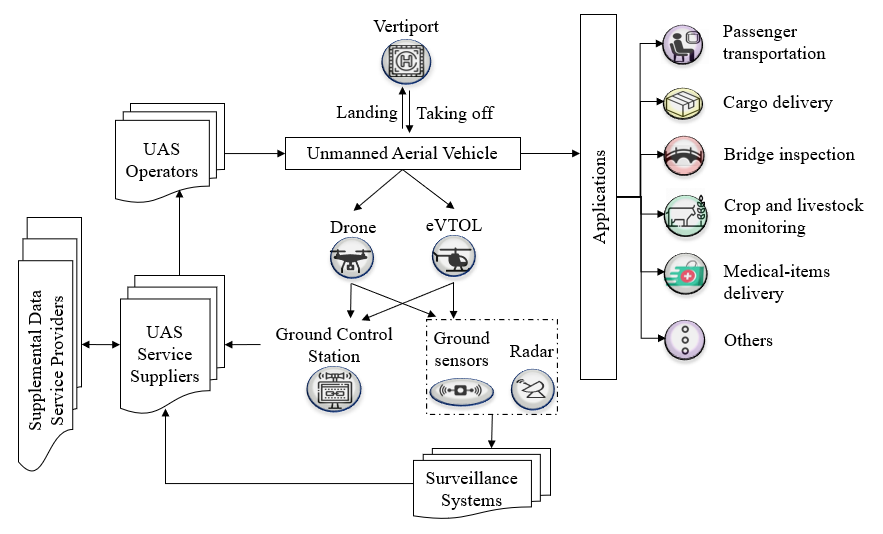}
    \caption{Notional AAM architecture.}
    \label{fig:my_label}
\end{figure}

An overview of the different primary stakeholders and components of the AAM architecture with their contextual linkages and functions is visualized in Figure \ref{fig:my_label}. We considered two types of unmanned aerial vehicle (UAV) or AAM aircraft in this architecture. The first one is the \textit{electric vertical takeoff and landing vehicle (eVTOL)
}, which have been developed to transport people and cargo through local, regional, inter-city, and urban spaces. eVTOLs are expected to significantly cut journey times relative to conventional ground transportation systems, while offering decreased pollution footprints, as well as increased operational safety, among other benefits \cite{Rothfeld2021}. The second one is the \textit{drone} or small uncrewed aircraft systems (sUASs), which has already been used to transport small packages quickly and effectively. In this paper, the terms 'drones' and 'sUASs' are used interchangeably. These vehicles have novel designs and capabilities, and are mostly intended to operate at altitudes less than 400 ft (120 m) above ground level in uncontrolled airspace \cite{kose2020simultaneous}. A \textit{vertiport} is considered to be an area designed for an UAV to perform take off, landing and recharging, which can be an airport, heliport, top of a building, or any other setup based on the location and type of application. To ensure safe UAS traffic management (UTM), the proposed AAM architecture is based on sharing of data through a distributed information network following a set of protocols and functions \cite{FAA2016}. The UAVs' telemetry data---such as position, velocity, flight intent, remote identification (RID), etc.---will be monitored by ground control stations via radio frequency link, and UAVs' movement will also be detected and tracked by the surveillance systems via ground sensors and radar. The data from ground control stations and surveillance systems will be routed to \textit{UAS Service Suppliers (USSs)}. The role of USS is to act as a third party UTM service provider for the UAS operators. They will offer information about the planned activities of airspace to the UAS operators so that they may plan and schedule their flights accordingly to accomplish their operations safely and efficiently. USSs will also be responsible for archiving operational data in historical databases for online and post-flight data analysis, regulatory, and operator accountability purposes \cite{FAA2020}. \textit{Supplemental Data Service Providers (SDSPs)} provide advanced flight support services, such as topographical and obstacle data, specialized meteorological data, surveillance, and constraint information to USSs. UAS operators will receive approvals, advisories and other UTM services from USSs and operate their fleet of UAVs accordingly for different applications, such as passenger and cargo transportation, bridge inspections, agriculture and livestock observation, medical delivery, etc. Thus, the physical and software infrastructure needed to enable these AAM applications include vertiports, ground control stations, surveillance systems, USSs and SDSPs. These components are collectively referred to as the \textit{AAM infrastructure}. 

The operational demands and projected advantages of AAM are pushing partnerships between public and commercial stakeholders, headed by the FAA and NASA, to develop AAM infrastructure to ensure the safe, secure, efficient, and equitable integration of UAS operations in the NAS \cite{FAA2016}. Additionally, a broad range of technical and regulatory challenges is currently being addressed by the government, academia, and industry as part of the AAM Mission to enable roll out of AAM in the near future \cite{AAMt}. The existing literature has mostly focused on concepts and technologies required for AAM operations; however, to successfully implement any revolutionary new technologically advanced transportation system, it is necessary at first to check if the implementation would be cost-effective. \textit{Cost-benefit analysis (CBA)} is an approach that validates the decision-making process for the government and private service providers and offers them a sound financial basis for their investments in such a system. CBA assists with determining if the payback or return on investment outweighs the costs and whether the project is financially robust and secure enough. It also streamlines the complicated investment decisions that must be made throughout the implementation. Before the government takes a step to invest in this new AAM infrastructure, the monetary and, especially, the non-monetary valuation of the AAM infrastructure should be taken under consideration to identify whether the investment has the potential to bring significant benefits and outcomes to the people of its state. To address this gap, we conducted a rigorous cost--benefit analysis of the AAM infrastructure using real world historical data, with a focus on identifying and evaluating the non-monetary benefits of AAM for the society and environment, from the perspective of the government for the 2022--2032 time period. The state of Ohio was chosen as the place for implementation of AAM infrastructure in this study. More specifically, we focus on the six major cities of Ohio---Akron, Cincinnati, Cleveland, Columbus, Toledo, and Dayton---as they have the highest market potential for AAM in Ohio. In this context, the following objectives were sought to be achieved by the authors of this paper:

1.	To identify and measure the most critical and substantial cost and benefit factors for society and environment associated with the AAM infrastructure implementation in Ohio in its six major cities.

2.	To determine if investment in AAM infrastructure is justified by evaluating its net positive gain.

The remainder of this paper is structured as follows: Section \ref{sec2} examines and explores the relevant literature. In Section \ref{sec3}, we discussed the approach and model used for forecasting various variables of interest for the next 11 years. The feasible benefit and cost factors related to AAM infrastructure have been identified and evaluated in \mbox{Sections \ref{sec4} and \ref{sec5}}, respectively.  After that, Section \ref{sec6} presents and discusses the results. Lastly, Section \ref{sec7} concludes the paper with the insights gained from the analysis and potential future work.

\section{Literature Review}\label{sec2}

This section is organized into two segments. The first segment reviews some key literature related to AAM infrastructure, operations, and regulations. Then, the next segment discusses the studies related to CBA in transportation projects for understanding the importance of CBA for implementing this emerging transportation infrastructure.

\subsection{Advanced Air Mobility Infrastructure, Operations, and Regulations}
Although UAS were previously mostly used by the military for intelligence, surveillance, reconnaissance, and combat missions \cite{Canis2015}, the demand for UAS for civilian use cases has risen dramatically in recent years. NASA's AAM vision is to assist the growing aviation markets in developing a safe air transportation system that carries people and freight using UAS currently being developed for this purpose. Substantial funds have been invested to build a ground infrastructure for AAM that includes ground sensors, vertiports, and the data exchange architecture for UTM system \cite{thipphavong2018urban}. During the 2020 fiscal year, the Aeronautics Research Mission Directorate (ARMD) established the AAM Mission Integration Office with the goal of promoting flexibility and agility as well as fostering collaboration across ARMD programs contributing to the AAM Mission \cite{Canis2015}. For a collaborative and responsible automation, the mission focuses on Automated Flight and Contingency Management (AFCM) subproject for developing Urban Air Mobility (UAM) system architectures and research findings to support standards for vehicle and pilot interface systems which will benefit the US industry \cite{AAMt1}. The National Campaign is a subproject under the AAM project designed to push the industry to perform ever more challenging ecosystem-wide system level safety and integration scenarios, show realistic and scalable system concepts, and establish a knowledge basis for requirements and standards \cite{AAMt2}. The High Density Vertiplex (HDV) subproject is in charge of developing and testing of concepts, requirements, software architectures, and technologies for the terminal environment surrounding vertiports, with a particular emphasis on how automation can improve the safety, efficiency, and scalability of flight operations in these environments \cite{AAMt3}.

For early users of civilian drone technology, there was once a legal gray area because laws and regulations governing the operation of these vehicles were essentially nonexistent, particularly at the federal and state levels \cite{Jiang2016}. Drone safety and efficiency were recognized as significant priorities for the near-term future when the FAA Modernization and Reform Act was passed in 2012. By 30 September 2015, the legislation demands the establishment of a strategy to integrate civil UAS into the NAS safely \cite{West2015}. NASA and the FAA have collaborated to set up a Research and Testing Task Force (RTT) to develop procedures that will allow large-scale operations of UAS \cite{FAA2016}. As seen, AAM infrastructure is a complex interaction between the FAA, the operator, and the different businesses providing services inside the ecosystem. It highlights a network that largely relies on the use of third-party entities \cite{FAA2020}. Because of the infrastructure, aircraft will no longer be required to communicate with a single entity, such as a designated air traffic controller. As an alternative, it will be able to communicate with a variety of service providers freely. 

Several businesses and academic researchers have been collaborating with NASA and the FAA to assist in integrating UAS into NAS. Airbus' department that focuses on developing digital traffic management infrastructures, Airbus UTM, is creating programs that will help establish the future AAM frameworks \cite{petrescu2017unmanned}. Amazon proposed assessing safe access with a best-equipped, best-served model for UAS, categorizing UAS based on vehicle systems and technologies, also known as equipage. Moreover, Amazon proposes a model based on a single operator controlling a wide range of vehicles, which substantially decentralizes decision-making authority among operators \cite{ali2019traffic}. Google's UAS Airspace System architecture divides the airspace into classifications as defined by the FAA \cite{gharibi2016internet}. Studies have also been done to build models of operations management, which are required to assist service providers in their operational decision-making, such as scheduling, dispatching, and fleet planning, in order to maximize their preferred objectives~\cite{shihab2019,Shihab2020}. According to the government's applicable regulatory, reliability, and performance requirements, these providers will be able to collaborate with the entire network to make cost-effective decisions based on flight missions. In addition, research has been conducted to maximize autonomous flight performance for tactical UAVs with flight control systems design \cite{oktay2017simultaneous}. 

\subsection{Cost--Benefit Analysis of Transportation Systems}
CBA is a quantitative and systematic procedure of evaluating a project's worth, which has become a preferred method in transportation projects for evaluating the impacts of travel time saving, safety increase, greenhouse gas reduction, improved air quality, etc.~\cite{Couture2016}. It assesses the effectiveness of transportation policies based on the monetary and non-monetary benefit factors. Policymakers have advanced a variety of arguments, and that is why they are optimistic about the incorporation of CBA, mainly in the planning phase regarding transportation projects \cite{Mishan2020}. According to NASA's market study, UAM has an economic viability for quicker transportation of people and packages \cite{hasan2018urban}. CBA was done for the Amazon Prime Air delivery system for utilizing drone technology to deliver items directly to consumers' doorsteps \cite{Welch2015}. This article examined how the benefits of labor-saving technology and speedier delivery outweigh the expenses of running the drone system. UAV offers a form of transportation that allows for moving air freight more cheaply and efficiently, and a CBA was performed to check its feasibility \cite{Groningen2017}. A methodology for migrating ATM operations to cloud computing to achieve cost savings and performance and efficiency benefits was developed by \cite{Ren2014}. An ideal but realistic assessment for cost-effective research could be used to comprehend the pluses and minuses of the upcoming air transportation system and to evaluate the installation of an air corridor within the Japanese airspace \cite{Barelli2018}. The project, presented in \cite{Adler2010}, seeks to provide a framework to analyze investments in high-speed rail and hub-and-spoke legacy airlines, including their impacts on transport equilibrium, and how the investments will impact the quality of life of the public. The merits of integrating electric automobiles into urban mobility were examined to determine its monetary viability \cite{Zito2004}. CBA can be used in transport planning to focus on travel time saving, which is essential in the evaluation of transportation projects because of their significance in reducing trip time \cite{Martens2017}. This analysis also can be done for air quality improvement reducing fuel consumption \cite{Bi2006}. Due to the rapid growth of the AAM and UAS industries, Ohio currently has lots of new possibilities and projects. Ohio Department of Transportation (ODOT) aims to implement AAM to ensure that people and cargo are transported safely, efficiently, and fairly across the state. A study was performed by the Crown team to analyze the economic impact of AAM in Ohio’s major urban centers~\cite{del2021infrastructure}. In this study, they analyzed some monetary benefit factors, such as GDP growth, tax revenue, passenger and cargo revenues.
 
After analyzing the literature, it can be said that CBA performs a significant role in the transportation field. Furthermore, the AAM infrastructure will be the next transportation system opening plenty of demands and opportunities. Though economic analyses were implemented in traditional transportation fields, it is not implemented for AAM with respect to non-monetary benefits like travel time savings, delivery delay reduction, reduction in carbon footprint, etc. Non-monetary benefit factors are equally important as the monetary factors but more difficult to evaluate. Typical studies applied CBA focusing on one factor, for example, UAS for safer bridge inspection quantifying its cost and benefit \cite{Hubbard2020}; but not covering the possible factors as much as possible in a study to get an overview of AAM's net positive gain. To address this research gap, we conducted a thorough data-driven quantitative cost–benefit analysis of AAM from the perspective of the state government to identify and estimate the most relevant and significant benefit and cost factors of AAM for society and environment.

\section{Forecasting Approach}\label{sec3}
\begin{figure}[H]
    \includegraphics[width=13cm,height=7cm]{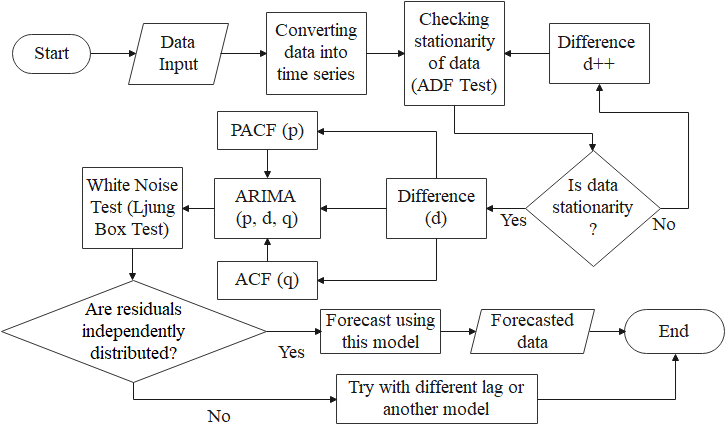}
    \caption{Logical flowchart for the ARIMA model.}
    \label{fig2}
\end{figure}

In CBA, forecasting is an integral tool for estimating the costs that will be occurred before incurring those. As this tool gives a possible preview of the future, it is commonly used in the transportation field during the cost--benefit analysis \cite{Martens2017}. To forecast the future values of 14 out of all variables (see Table \ref{Table1}) in our study, we performed a time-series analysis over 11 years period using an Auto-regressive Integrated Moving Average (ARIMA) model in RStudio. A time series is a sequence of data points that have been listed in time order. ARIMA model is considered as a powerful tool for doing forecasting in transportation~\cite{miller2019arima}. The `auto regressive' (AR) portion of ARIMA denotes that the dependent variable is regressed on its own lagged values. The middle part `integrated' (I) implies that the data have been substituted with the difference between their current and former values. The `moving average' (MA) section of the model denotes a linear regression of errors whose values existed simultaneously and at different periods in the past \cite{arima}. Figure \ref{fig2} shows the logical flow that we followed to build our ARIMA forecasting model. At first, we converted historical data of the corresponding variables into time series and then checked the stationarity of the time series. For stationarity, the time series must have constant mean, variance, and auto-correlation. Here, the auto-correlation of a time series measures the statistical relationship between its actual and lagged values. For checking the stationarity of the time series, we performed the Augmented Dickey--Fuller test (ADF test), where the null hypothesis ($H0$) is that the time series is not stationary; and the alternative hypothesis ($Ha$) is that the time series is stationary.

\vspace{-6pt}
\begin{table}[H]
\caption{Summary of forecasted variables using the ARIMA model.}
\label{Table1}
\tabcolsep=0.545cm
\begin{tabular}{cccc}
\toprule
\multicolumn{2}{c}{\textbf{Entity}} & \textbf{ARIMA(\emph{p}, \emph{d}, \emph{q})} & \textbf{Ljung-Box Test} \\ \midrule
\multicolumn{2}{c}{Vehicle Miles Traveled (\textit{\textbf{VMT
}})} & \multicolumn{1}{c}{ARIMA(0, 2, 1)} & \multicolumn{1}{c}{\emph{p}-value = 0.5683} \\ \midrule
\multicolumn{2}{c}{Population (\textit{\textbf{P}})} & \multicolumn{1}{c}{ARIMA(0, 2, 1)} & \multicolumn{1}{c}{\emph{p}-value = 0.9822} \\ \midrule
\multicolumn{2}{c}{Value of Statistical Life (\textit{\textbf{VSL}})} & \multicolumn{1}{c}{ARIMA(0, 2, 0)} & \multicolumn{1}{c}{\emph{p}-value = 0.8907} \\ \midrule
\multicolumn{2}{c}{Median Household Income (\textit{\textbf{MHI}})} & \multicolumn{1}{c}{ARIMA (0, 1, 0)} & \multicolumn{1}{c}{\emph{p}-value = 0.6957} \\ \midrule
\multicolumn{1}{c}{\multirow{3}{*}{Harvested area}} & \multicolumn{1}{c}{Soybean (\textit{\textbf{$A_S$}})} & \multicolumn{1}{c}{ARIMA(0, 0, 0)} & \multicolumn{1}{c}{\emph{p}-value = 0.9131} \\ \cmidrule{2-4} 
\multicolumn{1}{c}{} & \multicolumn{1}{c}{Corn (\textit{\textbf{$A_C$}})} & \multicolumn{1}{c}{ARIMA(0, 0, 0)} & \multicolumn{1}{c}{\emph{p}-value = 0.7899} \\ \cmidrule{2-4} 
\multicolumn{1}{c}{} & \multicolumn{1}{c}{Wheat (\textit{\textbf{$A_W$}})} & \multicolumn{1}{c}{ARIMA(0, 0, 0)} & \multicolumn{1}{c}{\emph{p}-value = 0.6069} \\ \midrule
\multicolumn{1}{c}{\multirow{3}{*}{Crop Yield}} & \multicolumn{1}{c}{Soybean (\textit{\textbf{$Y_S$}})} & \multicolumn{1}{c}{ARIMA(1, 1, 3)} & \multicolumn{1}{c}{\emph{p}-value = 0.3430} \\ \cmidrule{2-4} 
\multicolumn{1}{c}{} & \multicolumn{1}{c}{Corn (\textit{\textbf{$Y_C$}})} & \multicolumn{1}{c}{ARIMA(1, 1, 2)} & \multicolumn{1}{c}{\emph{p}-value = 0.2767} \\ \cmidrule{2-4} 
\multicolumn{1}{c}{} & \multicolumn{1}{c}{Wheat (\textit{\textbf{$Y_W$}})} & \multicolumn{1}{c}{ARIMA(3, 1, 1)} & \multicolumn{1}{c}{\emph{p}-value = 0.7521} \\ \midrule
\multicolumn{1}{c}{\multirow{3}{*}{Price}} & \multicolumn{1}{c}{Soybean (\textit{\textbf{$p_S$}})} & \multicolumn{1}{c}{ARIMA(0, 0, 0)} & \multicolumn{1}{c}{\emph{p}-value = 0.8212} \\ \cmidrule{2-4} 
\multicolumn{1}{c}{} & \multicolumn{1}{c}{Corn (\textit{\textbf{$p_C$}})} & \multicolumn{1}{c}{ARIMA(1, 0, 1)} & \multicolumn{1}{c}{\emph{p}-value = 0.3311} \\ \cmidrule{2-4} 
\multicolumn{1}{c}{} & \multicolumn{1}{c}{Wheat (\textit{\textbf{$p_W$}})} & \multicolumn{1}{c}{ARIMA(2, 1, 1)} & \multicolumn{1}{c}{\emph{p}-value = 0.0925} \\ \midrule
\multicolumn{2}{c}{Livestock count (\textit{\textbf{L}})} & \multicolumn{1}{c}{ARIMA(0, 1, 1)} & \multicolumn{1}{c}{\emph{p}-value = 0.3429} \\ \bottomrule
\end{tabular}
\end{table}

If $p$-value is less than or equal to 0.05 (choosing the confidence level as 95\%), we rejected the Null referring to the time series is stationary. 
Conversely, if the $p$-value is greater than 0.05, the time series is not stationary. Then, we did differencing
time series and repeat the ADF test until we obtain a stationary time series. We selected the final $d$ value. Next, PACF (Partial Auto-correlation Function) and ACF (Auto-correlation Function) were plotted, and the \emph{p}- and \emph{q}-values were selected accordingly. \emph{p} indicates the values of correlation between observations of a time series with their lagged values, whereas \emph{q} indicates the correlation between errors of a time series with their lagged values. The ARIMA(\emph{p}, \emph{d}, \emph{q}) model means that it can describe the forecasted data by combining a \emph{p} order auto-regressive model and a \emph{q} order moving average model, when the input data have gone through \emph{d} order differencing~\cite{stram1986temporal}. After that, we ran the ARIMA($p$, $d$, $q$) model and performed a white noise test with a Ljung Box Test for residuals. Here,

\begin{quote}
\textbf{Null Hypothesis} (\boldmath{$H$\textbf{0}})\textbf{:} The residuals are independently distributed. 
\end{quote}

\begin{quote}
\textbf{Alternative Hypothesis} (\boldmath{$Ha$})\textbf{:} The residuals are not independently distributed (they exhibit serial correlation).
\end{quote}

Generally, a model is considered as the best model if the residuals are independently distributed. Thus, if the $p$-value is greater than 0.05, it means this model is good for forecasting our time series. If not, different lag or another model has to be tried. Finally, we forecasted the time series with the fitted model. While forecasting each variable, we used a 95 percent confidence interval to consider uncertainty in those variables.

\section{Benefit Factors}\label{sec4}
In this study, nine benefit factors (BFs) were covered as part of the CBA, which highlight the savings in cost and time in various use cases that could be potentially achieved by AAM. We calculated the worth of the relevant non-monetary benefit factors in three steps. Identifying the benefit factors was the first step. Then, the positive impacts of these factors were measured in either physical or time-based units, and, finally, we translated the physical or time-based quantities into monetary value. The following sections describe the impact of each benefit factor identified as part of this study.

\subsection{BF-1: Travel Time Savings for Passengers}\label{sec 4.1}
This section examines the value of travel time, and travel time savings resulting from AAM passenger transportation using eVTOLs. Travel time is one of the most significant transportation costs, and travel time savings are often the primary justification for transportation infrastructure improvements \cite{litman2009transportation}. Various studies have developed estimates of travel time values for different user types and travel conditions. The Department of Transportation (DOT) and other government agencies always take initiatives that are intended to lower the travel time of the passengers. \textit{Value of Travel Time Savings (VTTS)} refers to the amount of money that travelers are willing to pay to reduce the travel time. In US DOT's last revised version of departmental guidance for value of travel time in economic analysis, the recommended VTTSs for local and intercity travel on surface mode (except high speed rails) were \$14.1 and \$20.4 per person-hour in 2015, respectively \cite{DOT:2016}. For our analysis, we considered the average of the two values (\$17.25 per person-hour) as VTTS. 

To estimate yearly VTTS for the forecast period, the factors on which VTTS depends were taken into account. According to the US DOT report, the valuation of travel time indicated that VTTS depends on the demographic characteristics of the traveling population (age, sex, hourly median income, personal interest, etc.), location (local or intercity), and the mode (surface, high-speed rails or air) and purpose of travel (personal or business). US DOT considered a set of weighted values for giving different priorities to different purposes of travel, which were derived from the 2001 National Household Travel Survey (NHTS). Considering surface modes, for local travel, the weighted values were 4.6\% for business travel and 95.4\% for personal travel; for intercity travel, the weighted values were 21.4\% for business travel and 78.6\% for personal travel. These same values were used in all the previous reports of DOT. Following their reports, it can be assumed based on the same city or state that all the variants mentioned above, except hourly median income, will remain more or less unchanged every year if there is no unusual situation. Using Equation (\ref{1}), we estimated the \textbf{VTTS' 
} which refers to the change in VTTS 
for the next 11 years relative to the year 2015 because \$17.25 was formulated based on the year 2015. For estimating the \textbf{VTTS'}, we used 31 years of historical data of yearly Median Household Income (MHI) of Six Major Cities in Ohio (SMCO) to forecast MHI for the next 11 years based on the process described in Section \ref{sec3}. Then, the forecasted values were placed in a column vector, \textbf{MHI}, with 11~$\times$~1 order. In this paper, all the bolded variables are the column vectors having 11~$\times$~1 order, whereas the non-bolded variables are regular variables. In Equation (\ref{1}), $i$ denotes the year index, $VTTS_i$ the VTTS for year $i$ and $MHI_i$ the MHI for year $i$:
\begin{equation} \label{1}
 VTTS'_i= \frac{MHI_i}{MHI_{2015}} \times VTTS_{2015}; \indent  i \in \{2022, 2023,..., n (year)\}
 \end{equation}

According to several studies, travel time-saving is significant if passengers travel by air using AAM for inter-regional trips, for example, city center to city center, or city center to rural proximity, rather than short inter-city trips \cite{del2021infrastructure}. The Crown team considered the travel trip-distance between the cities for regularly scheduled or on-demand transportation to be over 50–75 m. We assumed the average distance of a trip ($d_1$) to be 50 m for our analysis.  
The difference between eVTOL flight and ground vehicle commute times was determined using ArcGIS real-time traffic data in \cite{del2021infrastructure}. A 50 m trip by eVTOL is estimated to take 30 min, whereas a ground vehicle takes 80 min long for the same trip. Thus, the trip-time saving ($s$) is 50 min for each passenger in each trip using eVTOL. If $d_1$ is considered more than 50 m, $s$ is estimated to be more than 50 min according to real-time traffic data. We collected the data of yearly estimated AAM passenger traffic (\textbf{D}) from \cite{del2021infrastructure}, and obtained $S_{p_i}$, which represents total time saved by estimated AAM passengers in year $i$ (unit in hour):
\begin{equation} \label{120}
S_{p_i}= \frac{D_i \times s}{60}  \indent  i \in \{2022, 2023,..., n (year)\}.
 \end{equation}

Using Equation (\ref{3}), $VTTS^T_{(i)}$ was calculated which refers to the aggregated VTTS coming from year $i$.
 \begin{equation} \label{3}
VTTS^T_{(i)}= S_{p_i} \times VTTS'_i \indent  i \in \{2022, 2023,..., n (year)\}
 \end{equation}

\subsection{BF-2: Safety Cost Reductions for Passengers}\label{sec 4.2}
The increasing number of traffic accidents is one of the leading causes of disability and fatality, putting significant economic and social strains on the affected commuters \cite{8941429}. Traveling by eVTOL is expected to be safer than by car, and this benefit can help eVTOLs to gain a significant level of public acceptance as a ride-sharing alternative \cite{garrow1}. One of the most major transportation concerns is safety, and cost reductions in this area are generally used as the fundamental criterion for upgrading the infrastructure \cite{garrow1}.

\textit{Vehicle Miles Traveled (VMT)} refers to the total quantity of travel for all the roadway vehicles in a certain geographic region during a given time period. \textbf{$VMT_{US}$} is the forecasted VMT for USA (see Table \ref{Table1}), where \textbf{VMT'} refers to the forecasted VMT for SMCO. Furthermore, \textbf{$P_{US}$} is the estimated population of USA for the next 11 years, which was taken from~\cite{ewert2015us}, and \textbf{P} is the forecasted population of SMCO. We acquired the historical VMT of USA data of 59 years and SMCO population data of 71 years from \cite{VMT, Population}, respectively, to compute \textbf{$VMT_{US}$} and \textbf{P} using ARIMA (see Table \ref{Table1}), respectively. Then, these variables were used to determine \textbf{VMT'}:
\begin{equation} \label{121}
VMT'_i= \frac{VMT_{{US}_{(i)}}}{P_{{US}_{(i)}}} \times P_i \indent  i \in \{2022, 2023,..., 2032 (year)\}.
 \end{equation}

According to Uber's report, the number of fatalities per 100 million passenger miles in ground transportation is denoted by $A_g$, which equals to 0.6, whereas the number of fatalities per 100 million passenger miles in on-demand urban air transportation is denoted by $A_a$, which is 0.3. So, based on the fatality rate per passenger mile, Uber aims to set eVTOLS's safety level as twice as driving a car \cite{8941429}. $F_i$ can be calculated using Equation (\ref{5}), which refers to the additional number of fatality projected to occur in year $i$ in SMCO's ground transportation relative to the on-demand urban air transportation:
\begin{equation} \label{5}
F_i = \frac{VMT'_i}{10^8} \times (A_g-A_a)  \indent  i \in \{2022, 2023,..., 2032 (year)\}
 \end{equation}

The implementation of AAM infrastructure would allow $F_i$ fatalities to be prevented in year $i$ if all the people using ground transportation switch to the new urban air transportation, which is highly unlikely to occur in the next 11 years. \textbf{R} refers to the possible amount of reduction in the number of fatalities when AAM infrastructure is implemented in SMCO. According to NASA, eVTOL would carry 1 to 4 passengers per trip \cite{edwards2020evtol}. Therefore, the estimated number of passenger trips by eVTOL in year $i$ can be calculated which is denoted by :
\begin{equation} \label{122}
D^T_i = \frac{D_i}{4}    \indent  i \in \{2022, 2023,..., 2032 (year)\}.
\end{equation} 

The product of $D^T_i$ and $d_1$ (see Section \ref{sec 4.1}) gives the number of miles traveled by the passengers in SMCO in year $i$ using eVTOL. Then, by dividing this number by $VMT'_i$ (see Section \ref{sec 4.1}), we obtained the ratio of vehicles miles traveled by eVTOL and existing ground transportation in SMCO, which was used to compute the \textbf{R}:
\begin{equation} \label{123}
R_i = \frac{D_i \times d_1}{VMT'_i} \times F_i  \indent  i \in \{2022, 2023,..., 2032 (year)\}.
\end{equation}

\textit{Value of Statistical Life (VSL)} is a measure of how much an individual is willing to pay to lower a specific degree of risk or how much risk an individual is ready to tolerate by not paying for the risk reduction. In Equation (\ref{7}), the total cost reduction in year $i$ for safety benefit ($CR^s_i$) was calculated using the values of \textbf{VSL}, which was forecasted (see Table \ref{Table1}) based on the past data taken from US DOT \cite{VSL}:
\begin{equation} \label{7}
CR^s_i= R_i \times VSL_i \indent  i \in \{2022, 2023,..., 2032 (year)\}
 \end{equation}

\subsection{BF-3: Cost and Delivery Time Savings in Package Delivery by Drones}\label{sec 4.4}

\subsubsection{Cost Savings for Logistics Companies}\label{subsubsec4.4.1}
The drone delivery system would help a company to save money and increase its operational efficiency \cite{sudbury2016cost}. Organizations would be able to focus their present human resources more toward creative, inventive endeavors that can help them to extend their services even more with drone-delivery. The employees will have more time to concentrate on the company's day-to-day operations while also assuring that consumers receive a greater quality product. This technology would allow people to spend more time and money on ideas that might one day transform the world, rather than performing repetitive operations that lower productivity levels unnecessarily \cite{yoo2018drone}. 

The global drone package delivery market is projected to grow from \$343.30 million in 2019 to \$39,013 million in 2030 at a 53.8\% Compound Annual Growth Rate (CAGR) \cite{drone}. Using Equation (\ref{12}), the value of global drone package delivery market in year $i$ ($V_i$) was calculated (unit in \$ million), relative to $V_{2019}$, which denotes the value of the global drone package delivery market (unit in \$ million) in 2019:
\begin{equation} \label{12}
V_i=V_{2019} \times (1+CAGR)^t   \indent i \in \{2022, 2023,..., 2032 (year)\}; t \in \{3, 4, ..., 13\}.
 \end{equation}
 
In 2019, North America achieved \$237.70 million in the drone package delivery industry, and 89\% share of this market value was captured by the US \cite{drone}. Moreover, it is estimated that the percentage will be increasing in the future. Thus, the market size of the US for drone package delivery $U_{2019}$ is \$211.553 million. The value of the US drone package delivery market relative to the value of the global drone package delivery market (r) is:
\begin{equation} \label{124}
r' = \frac{U_{2019}}{V_{2019}} = \frac{211.553}{343.303} = 0.6162.
\end{equation}

Using $r'$, we calculated $U_i$, which is the size of the drone package delivery market in the US in year $i$:
\begin{equation} \label{126}
U_i=r' \times V_i \indent i \in \{2022, 2023,..., 2032 (year)\}.
\end{equation}

Amazon has already been testing Prime Air delivery drones for several years \cite{AmazonDrone}. In the United States, Amazon ships over 2.5 billion packages every year \cite{drone3}, and, by 2022, it will be around 6.5 billion parcels \cite{amazon}. Amazon says that 87\% of its parcels weigh less than five pounds, making them ideal for the type of drone delivery that the FAA has approved~\cite{sudbury2016cost}. Thus, the number of packages that drones per year can potentially ship by 2022 in the USA ($a$) is 5.59 billion (0.86 $\times$ 6.5 billion). Using Equation (\ref{15}), we calculated $N_i$ which denotes the normalized value of USA market size for drone package delivery in year $i$ (\$ million) relative to the value of $U$ in year 2022:
\begin{equation} \label{15}
N_i=\frac{U_i*r'}{U_{2022}} \indent i \in \{2022, 2023,..., 2032 (year)\}
 \end{equation}

Then, the number of packages to be delivered or the number of annual drone delivery trips (as one drone can deliver one package at a time) in SMCO in year $i$ ($p_{{SMCO}_i}$) was calculated using Equation (\ref{16}), where $r_{p_i}$ is the ratio of \textbf{P} and \textbf{$P_{US}$} (see Section \ref{sec 4.2}). As the value of $a$ was given for year 2022, we calculated the normalized factor $N_i$ to obtain \textbf{$p_{SMCO}$} (unit in billion):
\begin{equation} \label{16}
p_{{SMCO}_i}=N_{i} \times a \times r_{p_i}\indent i \in \{2022, 2023,..., 2032 (year)\}
 \end{equation}

For estimating the number of working days of drone services in a year, we considered that Ohio averages 131 days of rainfall per year \cite{rainy}, and on a day with precipitation, there is a 61\% chance of a thunderstorm or moderate rain \cite{sudbury2016cost}. Drones have trouble flying in moderate rain or thunderstorms; therefore, they are likely to face 81 (131 $\times$ 0.61) days of flight delays or cancellations per year, leaving 284 uninterrupted operational days in a year. Package delivery drones, such as Amazon Prime Air delivery drones, have a maximum speed of 50 m per hour or 80.47 km per hour and maximum flight time of 30 min with 5 lb (or 2.268 kg) payload \cite{sudbury2016cost}. We considered 7.5 m (or 12.07 km) as the average flight radius ($d_2$) for our analysis. Now, 50 m/ 1 h = 50 m/ 60 min = 7.5 m/ 9 min. As a result, this drone can deliver an item to 7.5 m distance and come back to station in 24 min (less than 30~min), which includes 6 min for flight acceleration and deceleration, package loading and unloading, and flight climb and descent. Thus, the maximum number of trips per year ($Tr'$) can be calculated Equation (\ref{17}), where $R_t$ denotes the round trip time (24 min):
\begin{equation} \label{17}
Tr'= \frac{60}{R_t} \times 24 \times 284
 \end{equation}
 
Now, $Y_i$ is the number of delivery drones estimated to be in operation to deliver $p_{{SMCO}_i}$ packages in SMCO in year $i$. At any given time, we assumed that another $Y_i$ drones will be recharging. Then, $Y_i$ can be calculated using Equation (\ref{18}):
\begin{equation} \label{18}
Y_i=\frac{p_{{SMCO}_i}}{Tr'} \indent  i \in \{2022, 2023,..., 2032 (year)\}.
\end{equation}
 
Each station will require $Z_i$ number of reserve drones in year $i$ (assuming approximately 25\% of drones in constant operation) to be ready on its busiest day, including backup for repairs or breakdowns to prevent running out of drones: 
\begin{equation} \label{180}
Z_i= 0.25 \times Y_i  \indent  i \in \{2022, 2023,..., 2032 (year)\}.
\end{equation}

Thus, total number of drones is
\begin{equation} \label{181}
X_i = 2Y_i + Z_i  \indent  i \in \{2022, 2023,..., 2032 (year)\}.
\end{equation}

For calculating the potential cost savings of logistics companies from using delivery drones, we compared the cost associated with the drone delivery system with the cost associated with the traditional delivery system using trucks, cars, or vans. We consider the following scenario: a logistics company like Amazon delivers its items to its customers through United Parcel Service (UPS). A UPS driver who delivers items for $m=10$ h a day earns \$25 per hour, and each day a driver can deliver up to $n=250$ items. The delivery cost by truck ($q$) is around \$30 per hour after factoring in petrol and tolls. It does not include Amazon's expenses of delivering items to UPS centers, gasoline, and other transportation costs, or any other contractual fees \cite{sudbury2016cost}. On the other hand, the capital cost of a drone used by Amazon is between \$1000 to \$3000 for the drone, \$2000 for software and equipment, and the annual operating cost $C_o$ (operator, maintenance, recharging, and other costs) is on average \$800 \cite{sudbury2016cost}. We considered the average capital cost $C_c$ is \$4000 for each drone. Thus, the cost savings $CS^l_i$ in year $i$ from the view of logistics companies is: 
\begin{equation} \label{20}
CS^l_i= [{{\frac{m \times q}{n}}-{\frac{(C_c+C_o) \times X_i}{p_{{SMCO}_i}}}}] \times {p_{{SMCO}_i}}   \indent  i \in \{2022, 2023,..., 2032 (year)\}.
\end{equation}

\subsubsection{Lead Time Savings}
A delivery drone's objective is to deliver packages quickly anywhere, especially in difficult-to-reach regions. It can be designed to transport packages from their base to a specified location, such as the address of the customer who requested the item. Based on 2018 industrial drones specifications, a drone's delivery time for a 5 m urban trip was calculated as 8 min to deliver package in front of a customer's doorstep, whereas the calculated delivery time by a car was 19 min \cite{drone1}. Based on this calculation, the drone delivery times are estimated to be 7 min in the year 2020 and 5 min in year 2025. According to \cite{drone1}, the estimated time savings are 12 min and 14 min, respectively. Thus, we considered the Average Time Savings (ATS) associated with each drone delivery trip as 13 min relative to ground delivery trips. Using the values of $p_{{SMCO}_i}$ (see Section \ref{subsubsec4.4.1}), $ATS$ and $VDTS$ (see Section \ref{subsubsec4.3.2}), we computed the Cost of Lead Time Savings ($CS^{ld}_i$) in year $i$: 
\begin{equation} \label{21}
CS^{ld}_i= p_{{SMCO}_i} \times \frac{ATS}{60 \times 24} \times VDTS \indent  i \in \{2022, 2023,..., 2032 (year)\}.
\end{equation}

\subsection{BF-4: Time and Inventory Cost Savings in Cargo Delivery by eVTOLs}\label{sec 4.3}

\subsubsection{Cost Savings in Inventory}\label{subsec 4.3.1}
Inventory warehousing cost is not a value-added activity as it raises the price of the goods or service without increasing the customer's value. The use of Just in Time (JIT) inventory system aids in the reduction of non-value-added inventory expenses \cite{goshime2019lean}. Manufacturers will benefit from the use of eVTOLs for cargo delivery because they will increase JIT deliveries, which will in turn minimize the warehousing expenses by reducing the warehouse inventory. By delivering high-value products using eVTOLs through logistics air corridors, distribution operators and logistics companies can supply freight and cargo to regional industries faster \cite{del2021infrastructure}. Suppose a logistics company can deliver $x$ number of products every year through the traditional approach. Using eVTOLs, it would be able to deliver $x$ number of products in less than a year. Let the time difference be $y$ months, which means $x$ number of products will not be stored for $y$ months in the warehouse. Thus, here, the cost reduction is the cost saving related to operating a warehouse \emph{y} months fewer. This cost saving can be used to store more products or support other activities related to the company.

We assumed $d_1$ miles as the average trip distance for cargo delivery (see Section \ref{sec 4.1}). From \cite{del2021infrastructure} report, we obtained the estimated cargo revenue for next 25 years. According to that report, carrying 1.2 million tons (1 US ton = 2000 lb) of cargo will produce more than \$2.1 billion in revenue over the next 25 years. In each trip, an eVTOL is expected to carry 50--1000 lb (average payload $p^e=525$ lb) commodities. From these data, we calculated the number of estimated AAM cargo trips (\textbf{T}) for the next 11 years for our analysis. Using $s$ from Section \ref{sec 4.1}, $S_{t_i}$, which refers to the delivery time saved by the annual AAM cargo trips (unit in month), can be calculated as:
 \begin{equation} \label{8}
 S_{t_i} = \frac{T_i \times s}{30 \times 24 \times 60} \indent  i \in \{2022, 2023,..., 2032 (year)\}.
 \end{equation}

\begin{table}[H]
\caption{Average rent and size of warehouses in SMCO}
\label{tab 2}
\tabcolsep=0.42cm
\begin{tabular}{ccc}
\toprule
\textbf{City} & \textbf{Average Base Rent (\$/(sf $\cdot$ month))} & \textbf{Average Size (sf)} \\ \midrule
Akron & \$0.90 & {48,813 
} \\ \midrule
Columbus & \$0.67 & 48,215 \\ \midrule
Cincinniati & \$0.77 & {42,879} \\ \midrule
Cleveland & \$0.82 & 42,571 \\ \midrule
Dayton & \$0.99 & {28,581} \\ \midrule
Toledo & \$0.61 & {26,726} \\ \midrule
Average   for all cities & r = \$0.79 & z = 39,631 \\ \bottomrule
\end{tabular}
\end{table}

Three metrics commonly used to calculate how much a warehouse will cost each month are: (i) monthly base rental charge, (ii) facility size (in square feet (sf)), and (iii) expected operational costs (also called NNN). Table \ref{tab 2} shows the average base rent $r$ (unit in \$) and average size $z$ (unit in sf) for a warehouse in SMCO \cite{ware1}. A \textit{Triple Net Lease (NNN)} is a leasing structure in which the tenant is liable for all running expenditures associated with a property, such as building insurance payments, utilities, HVAC (heating, ventilation, and air conditioning) maintenance, and property taxes \cite{ware}. We considered NNN for a warehouse ($o$) as \$0.25 per sf per month \cite{w4}, and the average yearly wage of a warehouse worker ($w_1$) as \$27,867 \cite{w3}. In 2017, the average area of workplace in America's warehouse dedicated to each employee ($a$) was 138 sf \cite{w2}. Using Equation (\ref{9}), the average cost to operate a warehouse ($W$) can be calculated:
 \begin{equation} \label{9}
 W=(r+o) \times z + \frac{z}{a} \times \frac{w_1}{12}
 \end{equation}
 
Using Equation (\ref{10}), we calculated $CS^c_i$, the cost savings in cargo delivery in year $i$ by using eVTOLs:
\begin{equation} \label{10}
 CS^c_i= S_{t_i} \times W  \indent  i \in \{2022, 2023,..., 2032 (year)\}
 \end{equation}

Truck-trailers are commonly used for transporting cargo within a state. The weight limit for these trucks operating on Ohio's interstate highways is 48,000 lb \cite{evtolcost}, with a payload capacity share relative to empty weight is 200\% \cite{payloadpercentage}. Thus, a full truck can carry around $p^t=$24,000 lb commodities per trip, whereas we considered the average payload of eVTOL as $p^e$. The average marginal cost per mile for truck-trailer in Ohio is $c^t=\$1.417$ including fuel cost, truck-trailer lease or purchase payment, repair and maintenance, truck insurance, tires, tolls, and driver wage \cite{hooper2018analysis}. The projected passenger price by NASA for eVTOL is \$6~\$11 per passenger per mile considering revenue factor and some cost drivers, such as pilot wage, maintenance cost (including labor), vertiport support or landing fees, battery and charging, aircraft and insurance, and other expenses \cite{evtolcost}. Henceforth, it is estimated that around $c^e= \$34$ will be the price per mile to deliver $p^e$ commodities by a eVTOL. Using Equation (\ref{1001}), we obtained the amount of extra cost incurred (\textbf{CI}) in eVTOL cargo delivery:
\begin{equation}\label{1001}
 CI_i= [\frac{c^t}{p^t} - \frac{c^e}{p^e}] \times p^e \times d_1 \times T_i  \indent  i \in \{2022, 2023,..., 2032 (year)\}
 \end{equation}

Using Equation (\ref{1003}), we computed the cost difference between the operational costs of truck-trailers and eVTOLs (\textbf{CD}):
\begin{equation}\label{1003}
CD_i =  CS^c_i - CI_i \indent  i \in \{2022, 2023,..., 2032 (year)\}
 \end{equation}
 
According to Equation (\ref{181}), $X$ is 2.25 times of $Y$ considering the recharging, weather and backup issues:
\vspace{12pt}
\begin{equation}\label{1002}
X_i= 2Y_i + 0.25Y_i = 2.25Y_i \indent  i \in \{2022, 2023,..., 2032 (year)\}
 \end{equation}
 
For these issues, it is highly unlikely to achieve $CD_i$ cost saving by logistic companies in year $i$. Therefore, we included a factor $f^e=\frac{1}{2.25}=0.44$ in Equation (\ref{1004}) for considering these issues for eVTOL as well. Then, we evaluated the total cost savings in cargo delivery by eVTOLs ($CS^T_i$) for year $i$:
\begin{equation}\label{1004}
CS^T_i = f^e \times CD_i   \indent  i \in \{2022, 2023,..., 2032 (year)\}
 \end{equation}

\subsubsection{Lead Time Savings}\label{subsubsec4.3.2}
The eVTOLs are envisioned to transport large goods and freight to and from airports, distribution hubs, and manufacturers, as well as end-users or consumers. There are five major entities in a supply chain network, namely, suppliers, manufacturers, distributors, retailers, and consumers or users. Once an order is placed by an entity, the time it takes for the entity to receive the goods is referred to as \textit{Lead Time}. Clearly, buyers desire shorter delivery lead times \cite{gawor2019customers}. Moreover, given the increasing demand for e-commerce, lead time has become an increasingly major quantity of interest due to customers' time constraints~\cite{10.1145/322796.322805}. Reduction in lead time from an integrated supply chain network indicates the use of faster shipping methods \cite{Read1995DiversificationBE}. To assess the monetary worth of the lead time, a survey-based study was performed, where the variables were based on the experience and characteristics of customers in the USA. Relative to a standard lead time of 3--4 days, the findings revealed that customers are prepared to spend an extra \$12.80 for delivery today, \$8 for delivery tomorrow, and \$2.78 for delivery in two days. On an aggregated average, delivery lead time was estimated to be worth \$3.61 per day \cite{Read1995DiversificationBE}. Based on this, we assumed that the aggregated average Value of Delivery Time Savings ($VDTS$) applicable for SMCO in the next 11 years is \$3.61 per day. Using the value of $S_{c_i}$ obtained from Equation (\ref{8}), we calculated the average Annual Value of Delivery Time Savings ($AVDTS_i$) for year $i$: 
\begin{equation} \label{11}
 AVDTS_i = (S_{c_i} \times 30) \times VDTS  \indent  i \in \{2022, 2023,..., 2032 (year)\}.
\end{equation}

\subsection{BF-5: Cost and Travel Time Delay Savings for Bridge Inspections}\label{sec 4.5}

\subsubsection{Travel Time Delay Savings for Passengers}
To maintain the structural integrity of bridges in Ohio, they are inspected regularly and on an emergency basis by ODOT. ODOT is increasingly operating sUASs for this purpose. Every year, the ODOT inspects 14,500 bridges at different levels of intensity. Snooper trucks with a bucket carrying inspectors close to the underside of a bridge are used for around 400 inspections every year \cite{del2021infrastructure}. We focused on these $I=400$ heavy inspections for our study. The use of snooper trucks mandates lane closures, and due to this reason, passengers have to face travel time delays. Value of travel time and travel time cost were discussed in Section \ref{sec 4.1}.  

According to the traffic flow analysis conducted by ODOT, the total average hourly volume per lane for Ohio State is 1200 passenger car per hour per lane ($\alpha$) is \cite{flow}. For our study, we considered one lane closure for bridge inspections. The traditional bridge inspection approach needs on average 8 h per inspection as Traditional Lane Closing Time ($TLCT$). In comparison, it takes only 4 h per inspection on average as Reduced Lane Closing Time ($RLCT$) if a drone is used. Using Equation (\ref{22}), the average number of yearly passengers facing lane closure ($P_1$) due to traditional bridge inspection in Ohio can be computed:
\begin{equation} \label{22}
P_1=I \times TLCT \times \alpha
\end{equation} 

Similarly, $P_2$ can be obtained which refers to the average number of yearly passengers facing lane closure due to bridge inspection by drones:
\vspace{12pt}
\begin{equation} \label{185}
P_2=I \times RLCT \times \alpha.
\end{equation}

Closing one lane limits the route capacity and causes traffic congestion by creating a bottleneck on the open lane. Based on a review of common practices across the country, the average delay value ($ADV$) was chosen to be 10 min per vehicle \cite{flow2}. Considering one passenger per vehicle, we obtained $D_T$ which refers to the yearly delay faced by $P_1$ using Equation (\ref{24}), and $D_D$ which refers to the yearly delay faced by $P_2$ using Equation (\ref{25}):
\begin{equation} \label{24}
D_T = P_1 \times ADV
\end{equation}
\begin{equation} \label{25}
D_D = P_2 \times ADV
\end{equation}

Lastly, we computed the total value of delay time savings in year $i$, $TS^D_i$, from Equation~(\ref{26}) using $VTTS_i$ (see Section \ref{sec 4.1}):
\begin{equation} \label{26}
TS^D_i = (D_T - D_D) \times VTTS_i \indent  i \in \{2022, 2023,..., 2032 (year)\}
\end{equation}

\subsubsection{Cost Savings in Bridge Inspections}
Implementation of drones in bridge inspections can result in significant cost savings. The infrastructure such as bridges and highways, public buildings, and locomotives are subjected to systematic and regular inspection. The most common type of checking process is based on visual inspection carried out by certified inspectors. It is time-consuming and costly, whereas drones are emerging as a potential alternative by ensuring that the activities for inspecting infrastructure will continue with less expensive operation \cite{Kalaitzakis2021}. 

Skyward researched to see the return on investment in Ohio from drone bridge inspections \cite{flow1}. Their analysis calculated that \$3143 is needed to inspect a bridge using only snooper considering core hours and \$4152 is needed considering nights or weekends. Core hours are the intervals of the day when all employees are required to work, and night or weekend means by the intervals when workers work at night or on the weekend, outside the core hours. On the other hand, \$522 (based on core hours) or \$735 (based on nights or weekends) is needed to inspect a bridge using drones (see Appendix \ref{appendix}). Employee-related costs such as health plan expenses, workman's compensation charges, and pension plan expenses are referred to as fringe benefits, which are considered as 45\% for both of the cases. ODOT estimated that 80\% inspections are done during core hours and 20\% during nights or weekends \cite{flow1}. We used their data for calculating the cost of 400 bridge inspections. The yearly cost for 400 bridge inspections using only snooper ($C_1$) is
\begin{equation} \label{262}
C_1 = 400 \times 80\% \times 3143 + 400 \times 20\% \times 4152.
\end{equation}

Though sUAS will not replace all the inspection cases, ODOT expects to handle nearly half of the 400 snooper truck inspections. Thus, out of 400, 200 bridge inspections can be done using drones, while the other 200 inspections will be carried out using snoopers. The yearly cost for 200 bridge inspections using drones ($C_2$) is 
\begin{equation} \label{260}
C_2 = 200 \times 80\% \times 522 + 200 \times 20\% \times 735.
\end{equation}

Using Equation (\ref{276}), we computed the yearly cost for 200 bridge inspections done by snooper trucks ($C_3$):
\begin{equation} \label{276}
C_3 = 200 \times 80\% \times 3143 + 200 \times 20\% \times 4152.
\end{equation}

Using Equation (\ref{30}), we determined the total yearly cost savings in bridge inspections ($CS^B$).
\begin{equation} \label{30}
CS^B = CR_1 - (CR_2 + CR_3)
\end{equation}

\subsection{BF-6: Cost Savings in Agriculture}\label{sec 4.6}
AAM aircraft can be used in agriculture (both crop and livestock farming) to save time and effort, which is one of the significant benefits of AAM. Grazing lands can be monitored using sUAS technology, which is a useful solution for monitoring livestock  \cite{Yuksel2020}.

\subsubsection{Increased Crop Yield}
Drones equipped with sensors and digital imaging capabilities provide farmers with a more detailed view of their fields, allowing them to maximize agricultural yields and farm efficiency \cite{ag4}. Drones have revolutionized the spraying applications to determine which parts of the field require more or less water, pesticide, and fertilizer. They also have advantages for assessing storm damage, crop monitoring by assessing health, density, and also locating and possibly treating a disease before it spreads. Drones are boosting precision agriculture to new heights, and combining with other technology, farmers can gain the most out of their crops, regardless of the resources \cite{delavarpour2021technical}. For evaluating the value of increased crop yield or production and cost savings, we found soybeans as the most grown crop in Ohio. In addition,  corn and wheat, grown for both animal feed and human consumption, were found in the top ten agriculture in Ohio \cite{ag5}. Moreover, these are the most widely planted crops in the USA. Thus, we kept our focus on these crops for our analysis. 

Although many farmers in the United States were among the first to exhibit interest in the possibilities of sUAS, the expenses and time required to master the systems proved to be larger than expected because many farmers were still suffering from the financial impacts of the 2008 financial crisis when the technology began to take hold \cite{del2021infrastructure}. For this reason, we calculated a factor to compensate for the savings generated from this benefit factor. According to Teal's estimation, only farmers who earn more than \$100,000 a year are likely to employ this technology \cite{del2021infrastructure}. According to the United States Department of Agriculture (USDA), the number of farms in Ohio is 77,805, and the number of farms that make \$100,000 or more is 13,893 \cite{ag1}. The Factor expressing the ratio ($F_{ag}$) of farmers who will take on the technology is 13,893/77,805 or 0.18.

From Ohio Agricultural Statistics \cite{ag3}, we collected the historical crop data for 69~years. We used the harvested area (acre), yield (bushel), and price (dollar) data of soybean, corn, and wheat to forecast data for the next 11 years using the ARIMA model (see Section~\ref{3}), and put the forecasted data in the corresponding vectors (see Table \ref{Table1}). Wheat yields would increase by $x=3.3\%$ for using drones, while soybean and corn yields would rise by $y= 2.5\%$ \cite{ag5}. Using Equation (\ref{31}), we estimated the annual Increased Production of Soybean (\textbf{$IP_S$}):
\begin{equation} \label{31}
IP_{S_{(i)}} = (Y_{S_{(i)}} + Y_{S_{(i)}} \times y\%) \times A_{S_{(i)}} \indent  i \in \{2022, 2023,..., 2032 (year)\}
\end{equation} 

Then, we obtained the yearly estimated Value of Increased Production for Soybean (\textbf{$V_S$}) using Equation (\ref{32}):
\begin{equation} \label{32}
V_{S_{(i)}} = IP_{S_{(i)}} \times {p_{S_{(i)}}} \times F_{ag} \indent  i \in \{2022, 2023,..., 2032 (year)\}
\end{equation}

Similarly, we calculated estimated Value of Increased Production for Corn ($V_{C_{(i)}}$) and Value of Increased Production for Wheat $V_{W_{(i)}}$) for year $i$. Using Equation (\ref{33}), we computed the total estimated Value of Increased Production ($VIP^T_i$) in year $i$ from soybean, corn and wheat fields:
\begin{equation} \label{33}
VIP^T_i = V_{S_{(i)}} +V_{C_{(i)}}+V_{W_{(i)}} \indent  i \in \{2022, 2023,..., 2032 (year)\}
\end{equation}

\subsubsection{Cost Savings in Crop Farming}
Drones can lower the labor and equipment cost significantly. According to Informa Economics and Measure, a drone service provider firm, corn farmers stand to benefit the most from aerial devices, with average savings of \$11.58 per acre denoting $d$, whereas the farmers of wheat and soybean would save \$2.57 and \$2.28 per acre, respectively \cite{ag5}. Here, the estimated cost savings $C_C$ by drones in a corn field in year $i$ is 
\begin{equation} \label{330}
C_{C_i}= A_{S_{(i)}} \times d \times F_{ag} \indent  i \in \{2022, 2023,..., 2032 (year)\}.
\end{equation} 

Similarly, we calculated the yearly cost savings for soybean ($C_{S_{(i)}}$) and wheat ($C_{W_{(i)}}$). Using Equation (\ref{35}), we computed $CS^{crop}_i$ the total estimated value of cost savings in year $i$ from soybean, corn, and wheat fields in Ohio:
\begin{equation} \label{35}
CS^{crop}_i = C_{S_{(i)}} +C_{C_{(i)}}+C_{W_{(i)}} \indent  i \in \{2022, 2023,..., 2032 (year)\}
\end{equation}

\subsubsection{Cost Savings in Livestock Monitoring}
Drones for livestock monitoring are picking up steam in many nations. Drones are employed to keep an eye on the herd, requiring less workforce and allowing the farmers to follow the animals better than ever before \cite{al2020drones}. Farmers have also realized that using this technology will save them cash, time, and energy. A case study was conducted to monitor the impact of sUAS on farmland of $l_c=980$ cattle and ewes \cite{ag2}. In that study, it was found that the duty of inspecting the hill sheep was decreased from 40 min to 10 min, and, eventually, it saved $ST_1=26$ man-hours 
 per year. For checking some other cattle and sheep grazing on grassland, the time lowered from 1.5 h to 30 min, which saved a total of $ST_1=104$ man-hours per year. Traditionally, workers used to monitor those animals either walking on foot or driving quad or truck, but drones saved a lot of labor time. The average farm labor in Ohio ($FL$) is on average \$13.93 per hour \cite{range}. According to USDA, the animal, requiring large farmland to graze, are cattle, sheep, and lamb. These animals were on our consideration during the study as focusing on larger farmland could be more effective than the other animals' farmland in terms of cost savings. Using Equation (\ref{36}), we calculated $E$, which indicates the estimated labor cost savings per year per animal:
\begin{equation} \label{36}
E = \frac{(ST_1+ST_2) \times FL}{l_c}
\end{equation}

Using the forecasted number of animal in a year Livestock Count ($L$) (see Table \ref{Table1}), we obtained the annual estimated Value of Labor Cost Savings (\textbf{VLCS}) using Equation (\ref{37}):
\begin{equation} \label{37}
VLCS_i=E \times L_i \times F_{ag} \indent  i \in \{2022, 2023,..., 2032 (year)\}
\end{equation}

\subsection{BF-7: Life-Saving Value of Drones in Medical Delivery for Patients}\label{sec 4.7}
Medical delivery differs significantly from regular package delivery due to the delicate nature of items, namely, blood, lab samples, organs, medications, etc. As sUAS provides faster transmission of test samples, patients can receive treatment more quickly \cite{del2021infrastructure}. Here, the drones are considered suitable for organ delivery when the necessary configurations are applied. All medical deliveries are important without any doubt; however, emergency cases like Out of Hospital Cardiac Arrest (OHCA) should be given as the highest priority while designing the network of drones used for medical deliveries. When OHCA sufferers receive early defibrillation, they have a better chance of survival. In this paper, we focused on Automatic External Defibrillators (AEDs) delivery by drones in response to its critical scenario in the USA. According to the Ohio Department of Public Safety (ODPS), in the United States, Sudden Cardiac Arrest (SCA) is the leading cause of mortality \cite{med2}. They mentioned that OHCA accounts for around 250,000 deaths, and the yearly incidence of SCA in North America is 55 per 100,000 people. Using \textbf{P} (see Section \ref{sec 4.1} and Table \ref{Table1}), we calculated the number of OHCA sufferers ($U_{{OHCA}_{(i)}}$) in year $i$ in SMCO:
\begin{equation} \label{439}
U_{{OHCA}_{(i)}}= \frac{55}{\text{100,000}} \times P_i \indent i \in \{2022, 2023,..., 2032 (year)\}.
\end{equation}

An optimized approach was done to deliver AED assisted drones to the OHCA patients, where the authors analyzed their results considering six cases of drone docking stations needed to build a drone-delivery network through a state \cite{bogle2019case}. In our report, we used their data of number of drone docking stations, and put those in \textbf{DSN}: 
\begin{equation} \label{438}
DSN =  [0, 50, 200, 500, 750, 1015]. 
\end{equation}

For the six cases of drone docking stations, we also took the data of percentage of survivors with respect to OHCA incidents (\textbf{$p_s$}) and cost per additional survivor relative to `no drone' case (\textbf{CAS}): 
\begin{equation} \label{500}
p_s =  [12.3\%, 12.9\%, 13.3\%, 13.8\%, 14\%, 14.4\%]
\end{equation}  
\begin{equation} \label{520}
 CAS = [0, \$14752, \$31905, \$55792, \$73160, \$76495].
 \end{equation} 

When the number of drones is zero, the percentage of survivors is the lowest. With the increase of drone docking stations through a state, the delivery time will be lower, resulting in more chances of saving more patients' lives. On the other hand, with the increase of drone docking stations, the cost of preparing the stations will also be increasing. In Equation (\ref{39}), $N_{S_{(ij)}}$ is a matrix of 11~$\times$~6 order denoting the number of survivors in the next 11 years for the six drone docking cases:
\begin{equation} \label{39}
N_{S_{(ij)}} = U_{{OHCA}_{(i)}} \times p_{s_{(j)}} \indent  i \in \{2022,..., 2032 (year)\};  j \in \{1, 2, ..., 6\}
\end{equation}

Then, we computed $N_{S_{(ik)}}$ using Equation (\ref{40}) which stands for the number of increased survivors relative to the case of zero drones (`no drone' case or 1st column of $N_{S_{(ij)}}$). Here, the order of $N_{S_{(ik)}}$ is 11~$\times$~5:
\begin{equation} \label{40}
N_{S_{(ik)}} = N_{S_{(ij)}} -N_{S_{(i1)}} \indent  i \in \{2022,..., 2032 (year)\};  j \in \{2, ..., 6\};  k = j-1
\end{equation}

$CSS$ represents the cost savings per survivors that we evaluated using $VSL_i$ (see Section \ref{sec 4.2} and Table \ref{Table1}) and $CAS_j$:
\begin{equation} \label{409}
CSS_{(ik)} = VSL_i - CAS_j \indent  i \in \{2022,.., 2032 (year)\};  j \in \{2, ..., 6\}; k = j-1.
\end{equation}

$TVSL$ is a matrix of 11~$\times$~5 order denoting the total value of saving statistical life calculated using Equation (\ref{42}):  
\begin{equation} \label{42}
TVSL_{(ik)}=CSS_{{(ik)} \times N_{S_{(ik)}}}  \indent  i \in \{2022,..., 2032 (year)\}; k \in \{1, 2,.., 5\}
\end{equation}

\subsection{BF-8: Increased Tax Income for Government}\label{sec 4.8}
The state of Ohio is well-positioned to expand and continue lucrative AAM operations in urban, suburban, and rural regions over the next 11 years. The government will receive a significant level of tax from the local, state, and federal levels, which can lead to increased government expenditure, allowing the state to invest more in infrastructure, economic and social programs, and other sectors. We acquired the estimated tax income generated by AAM from \cite{del2021infrastructure}. In this paper, it is the only monetary benefit factor we considered. The government itself will directly receive this benefit, the expected amount of which serves as a strong potential reason to invest in this infrastructure.

\subsection{BF-9: Savings in Social Cost of Greenhouse Gases}\label{sec 4.9}
Scientists and the general public have been concerned about the possible repercussions of climate change for many years now. Conventional fossil Fuel-powered Ground Vehicles (CFGV) are one of the most significant sources of carbon emissions as they burn tons of gallons of fuel every year, resulting in global warming, increasing pollution levels and a variety of consequent health issues. Switching to eVTOLs from CFGVs can assist in alleviating the problem of excessive carbon emissions. It offers a new paradigm shift in green transportation, striving for a carbon-neutral and sustainable environment. \textit{Social Cost of Carbon (SCC)}, \textit{Social Cost of Methane (SCM)}, and \textit{Social Cost of Nitrous Oxide (SCN)} are used for monetary assessment of the social damage caused by emitting one metric ton of carbon dioxide ($CO_2$), methane ($CH_4$) and nitrous oxide ($N_2O$), respectively. These values reveal how much it is worth to us to avert future damage to our environment by avoiding the emission of these gases. After inflation adjustments, the interim $SCC$ was valued at \$51 per metric ton $CO_2$ by the Biden administration using a 3\% average discount rate ($r$) in 2020 \cite{e4}. The 2020 global $SCM$ was estimated to be \$1200 per metric ton $CH_4$ using a 3\% discount rate by the inter-agency working group on \textit{Social Cost of Greenhouse Gases (SCGHG)}, United States Government; and, for $SCN$, the value was estimated to be \$1500 per metric ton $N_2O$ \cite{e1}. Using 3\% discount rate, we forecasted the future value of SCC (\textbf{$F_{SCC}$)} for year $i$: 
\begin{equation} \label{43}
 F_{{SCC}_{(i)}} = SCC \times (1+r)^t \indent  i \in \{2022,..., 2032 (year)\}; t \in \{2, 3, ..., 12\}.
\end{equation}

Similarly, we obtained \textbf{$F_{SCM}$} for $CH_4$ and \textbf{$F_{SCN}$} for $N_2O$ using Equations (\ref{44}) and (\ref{45}), respectively:
\begin{equation} \label{44}
 F_{{SCM}_{(i)}} = SCM \times (1+r)^t \indent  i \in \{2022,..., 2032 (year)\}; t \in \{2, 3, ..., 12\}
\end{equation}
\begin{equation} \label{45}
 F_{{SCN}_{(i)}} = SCN \times (1+r)^t \indent  i \in \{2022,..., 2032 (year)\}; t \in \{2, 3, ..., 12\}
\end{equation}

The average social cost of \textbf{$F_{SCM}$} and \textbf{$F_{SCN}$} is \textbf{$F_{SCMN}$}:
\begin{equation} \label{46}
 FSC_{MN_{(i)}}= \frac{F_{{SCM}_{(i)}} + F_{{SCN}_{(i)}}}{2}  \indent  i \in \{2022,..., 2032 (year)\}.
\end{equation}

In 2018, the weighted average combined fuel mileage ($MPG$) of vehicles and light trucks was 22.5 m per gallon \cite{e2}. Then, using \textbf{VMT'} (see Section \ref{sec 4.2}) and $MPG$, we measured $N_{g(i)}$, which refers to the estimated number of gallons of fuel burned in year $i$ in SMCO:
\begin{equation} \label{600}
N_{g_{(i)}} = \frac{VMT'_i}{MPG} \indent  i \in \{2022,..., 2032 (year)\}.
\end{equation}

Carbon dioxide emissions as a proportion of total GHG emissions---including carbon dioxide, methane and nitrous oxide, emissions---for passenger vehicles ($F_{eq}$) was $0.993$ in 2018. The amount of $CO_2$ produced from burning each gallon of gasoline ($g_C$) is $8.89$~$\times$~$10^{-3}$ metric tons \cite{e2}; and the amount of $CH_4$ and $N_2O$ produced for burning each gallon of gasoline ($g_{MN}$):  
\begin{equation} \label{48}
 g_{MN} = g_C  \times  (\frac{1}{F_{eq}}-1).
\end{equation}

Thus, the amount of $CO_2$ (in metric ton) emitted in year $i$ in SMCO (\textbf{$A_{C}$}) is: 
\begin{equation} \label{49}
 A_{C_{(i)}} = g_C  \times  N_{g_{(i)}} \indent  i \in \{2022,..., 2032 (year)\}.
\end{equation}

Similarly, the amount of combined $CH_4$ and $N_2O$ (in metric ton) emitted in year $i$ in SMCO (\textbf{$A_{MN}$}) is: 
\begin{equation} \label{654} 
A_{MN_(i)}=g_{MN} \times  N_{g_{(i)}} \indent  i \in \{2022,..., 2032 (year)\}.
\end{equation}

The use of eVTOLs instead of CFGV for transportation will help to reduce the amount of emissions generated from the existing ground transportation system because eVTOLs have lower life cycle emissions than their ground counterparts. There are five stages in the life cycle of an eVTOL: (i) the manufacturing of components of the vehicle, (ii) the manufacturing of components of batteries, (iii) transportation of the vehicle from the assembly plant to the user, (iv) vehicle use, and (v) end of life \cite{andre2019robust}. Emissions from a CFGV is mostly produced during the fourth stage, the vehicle use stage; on the other hand, emissions from a eVTOL will mostly be caused during its manufacturing, the first and second stages. During the vehicle use stage of an eVTOL, they are estimated to have zero tailpipe emissions but a small amount of emissions will be associated with the recharging of the eVTOL, which will be produced at the power plant for electricity generation. Relative to CFGV, an eVTOL are envisioned to emit 35\% lower GHG, which will be more environment friendly \cite{kasliwal2019role}.

The total savings in social cost of GHG emissions potentially achieved by enabling AAM depends on (i) the amount of GHG emissions of eVTOLs and (ii) the number of trips travelled by eVTOLs instead of CFGVs. To take this into account, we considered two factors, fuel Factor ($F_f$) and demand factor $D_{f_{(i)}}$. Based on the above estimates, we considered $F_f=0.35$ denoting the ratio of reduced GHG emission by an eVTOL with respect to a CFGV.

Americans take around $Tr_{US}=411$ billion trips in a year using cars or station wagons, vans, and light trucks \cite{americantrips}. Total number of trips in SMCO using these vehicles ($Tr_{SMCO}$) in year $i$ can be computed using Equation (\ref{5300}):
\begin{equation} \label{5300}
 Tr_{{SMCO}_{(i)}} = \frac{P_i}{P_{{US}_{(i)}}} \times Tr_{US}  \indent  i \in \{2022,..., 2032 (year)\}.
\end{equation}

Recall \textbf{$D^T$}, \textbf{$p_{SMCO}$}, and \textbf{T} from Equations (\ref{122}) and (\ref{16}) and Section \ref{sec 4.4}, respectively.

Then, we computed the demand factor $D_{f_{(i)}}$, the ratio of estimated total number of trips travelled by eVTOL and ground vehicles in SMCO:
\begin{equation} \label{55}
 D_{f_{(i)}} = \frac{D^T_i+p_{{SMCO}_{(i)}}+T_i}{Tr_{{SMCO}_{(i)}}}.  \indent  i \in {2022,..., 2032 (year)}
\end{equation}

Lastly, $TSCS_i$ denotes the total social cost of GHG savings in year $i$:

\begin{adjustwidth}{0cm}{0cm} \begin{equation} \label{61}
TSCS_{(i)} = D_{f_{(i)}}  \times F_f \times (F_{{SCC}_{(i)}} \times A_{C_{(i)}}+FSC_{MN_{(i)}} \times A_{MN_{(i)}}) \indent  i \in \{2022,..., 2032 (year)\}.
\end{equation} \end{adjustwidth}

\section{Cost Factors}\label{sec5}
Our cost analysis is based on the findings of the CCO report \cite{del2021infrastructure}. This study has utilized over 100 sources of information in estimating capital investments and operating expenditures. They, first, categorized the cost factors into three groups: (i) Ground infrastructure, 
(ii) Systems for Providers of Services for Unmanned Air Vehicles (PSU), and (iii) Vehicles. 

Then, they estimated the corresponding Capital Expenditure (CAPEX) and Operating Expenditure (OPEX) for the first two categories and the acquisition and maintenance cost for the third category. CAPEX will be used for the feasibility study of the planning, ground site, environment analysis, power grid, communication networks, and others. 

The initial CAPEX is needed for land acquisition, site construction, building infrastructure, passenger amenities, safety facilities, and many more. In addition, upgrading existing facilities (e.g., helipads, heliports, and airports) to support eVTOLs' landing, refueling, recharging, and parking is also considered as a part of CAPEX. Some PSU facilities go under the CAPEX section; for example, purchasing different kinds of sensors, navigation beacons, radars, computers, and other types of equipment are included as a part of CAPEX.

OPEX consists of various recurring costs such as managing facilities, maintaining communication networks, labor costs, power bills, rental costs, and many. This cost section also includes subscription charges for Big Data analytics software, automation system, and flight decision support system. On top of that, the annual salary for computer operators, database managers, specialists in ATC and data analysis, IT technicians, and other professional personnel are included as OPEX cost. 

\begin{figure}[H]
    \includegraphics[width=13cm,height=6cm]{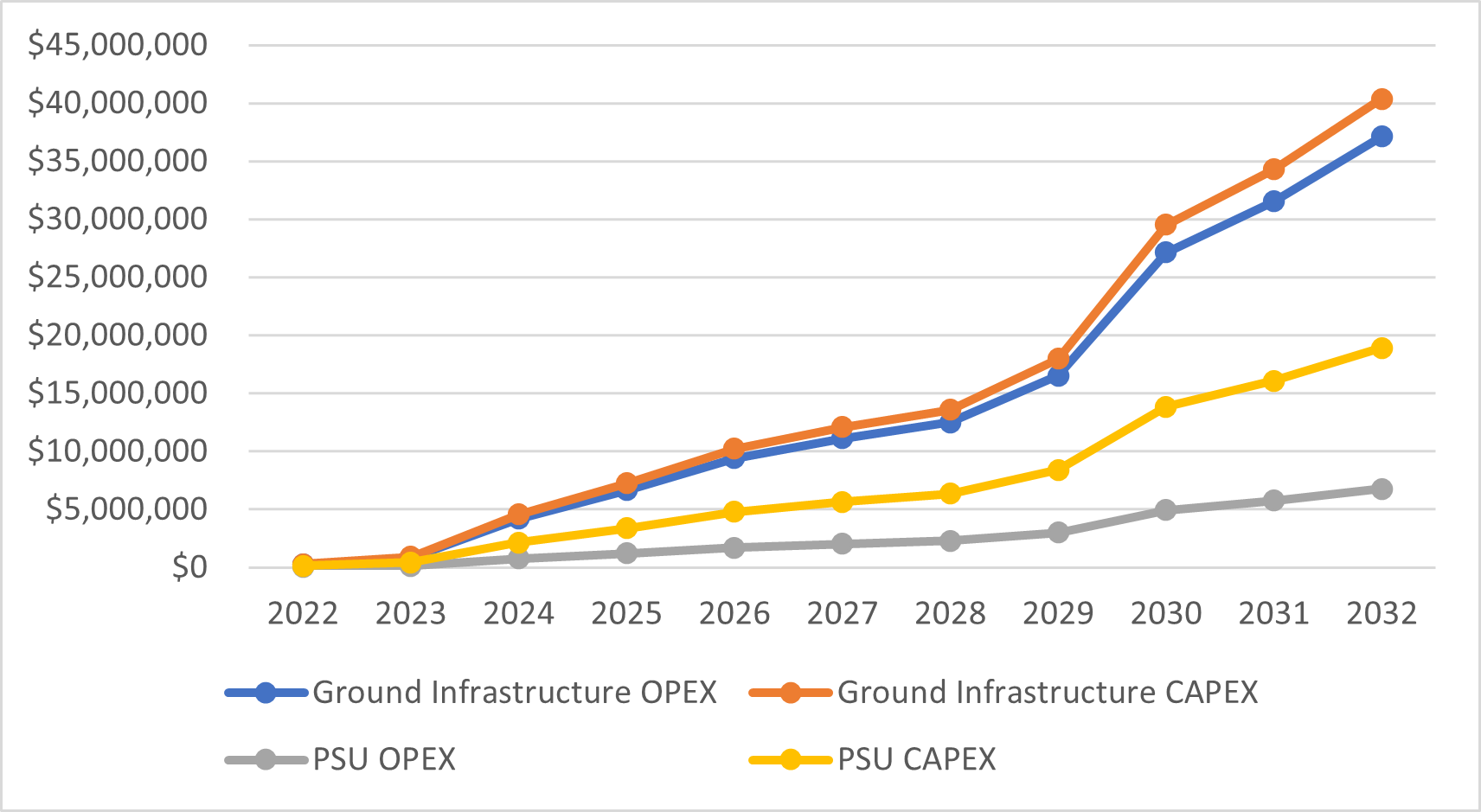}
    \caption{CAPEX-OPEX 
 costs for AAM infrastructure in SMCO (data taken from \cite{del2021infrastructure}).}
    \label{fig3}
\end{figure}

As we are conducting our analysis from the perspective of the state government of Ohio, we excluded the costs associated with the vehicle category since these costs will be incurred by the AAM operators, not by the government. Thus, in our analysis, we included the first two types of costs, which are presented in Figure \ref{fig3}. Though the CAPEX and OPEX costs could be funded by state initiatives, private sector investment, or a combination through public-private partnership models, we analyzed the worst-case scenario in our study, assuming the government will fund all these costs. If other entities fund these costs in a partnership with the government, the profit will eventually go up.

\section{Results}\label{sec6}
Using all the data and equations described in Sections \ref{sec4} and \ref{sec5}, we computed the dollar value of the benefits and costs associated with AAM infrastructure over the 2022--2032 time period. The estimated yearly net positive gain derived from AAM has also been determined. These results are presented in this section. Note that all the benefit factors are not illustrated in the same figure as their scales vary widely.

\begin{figure}[H]
    \includegraphics[width=13cm,height=6cm]{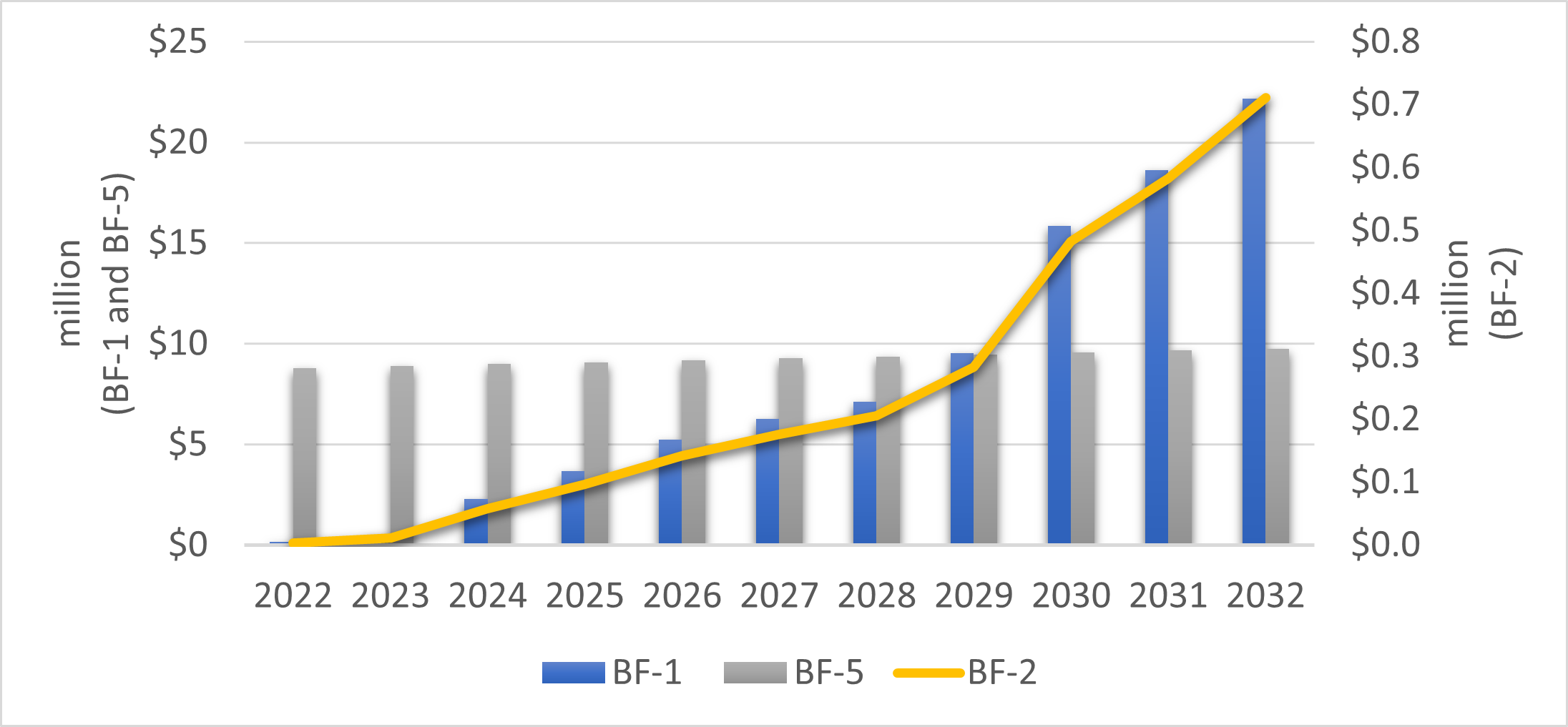}
    \caption{Yearly estimated values coming from BF-1 (see Section \ref{sec 4.1}), BF-2 (see Section \ref{sec 4.2}) and BF-5 (see Section \ref{sec 4.5}).}
    \label{fig4}
\end{figure}

\begin{figure}[H]
    \includegraphics[width=13cm,height=6cm]{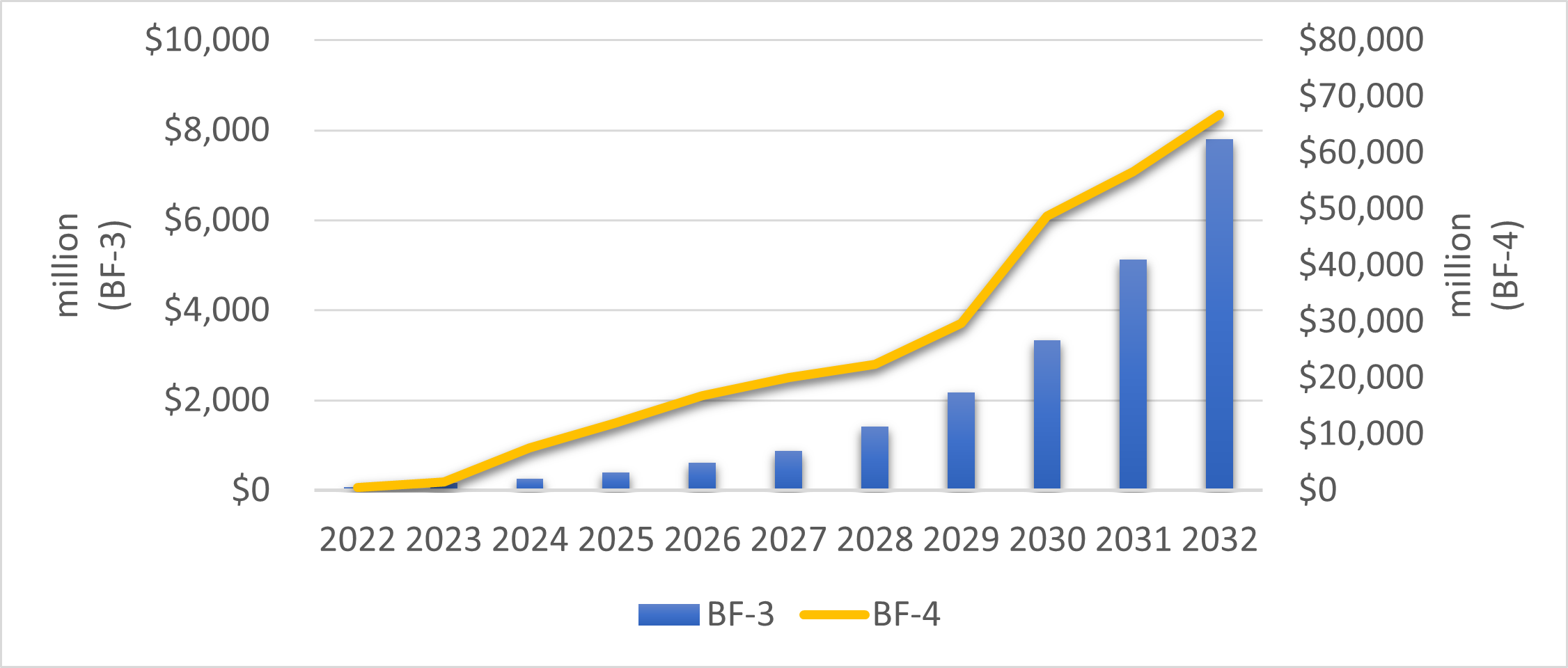}
    \caption{Yearly estimated values coming from BF-3 (see Section \ref{sec 4.4}) and BF-4 (see Section \ref{sec 4.3}).}
    \label{fig5}
\end{figure}

\begin{figure}[H]
    \includegraphics[width=12cm,height=6cm]{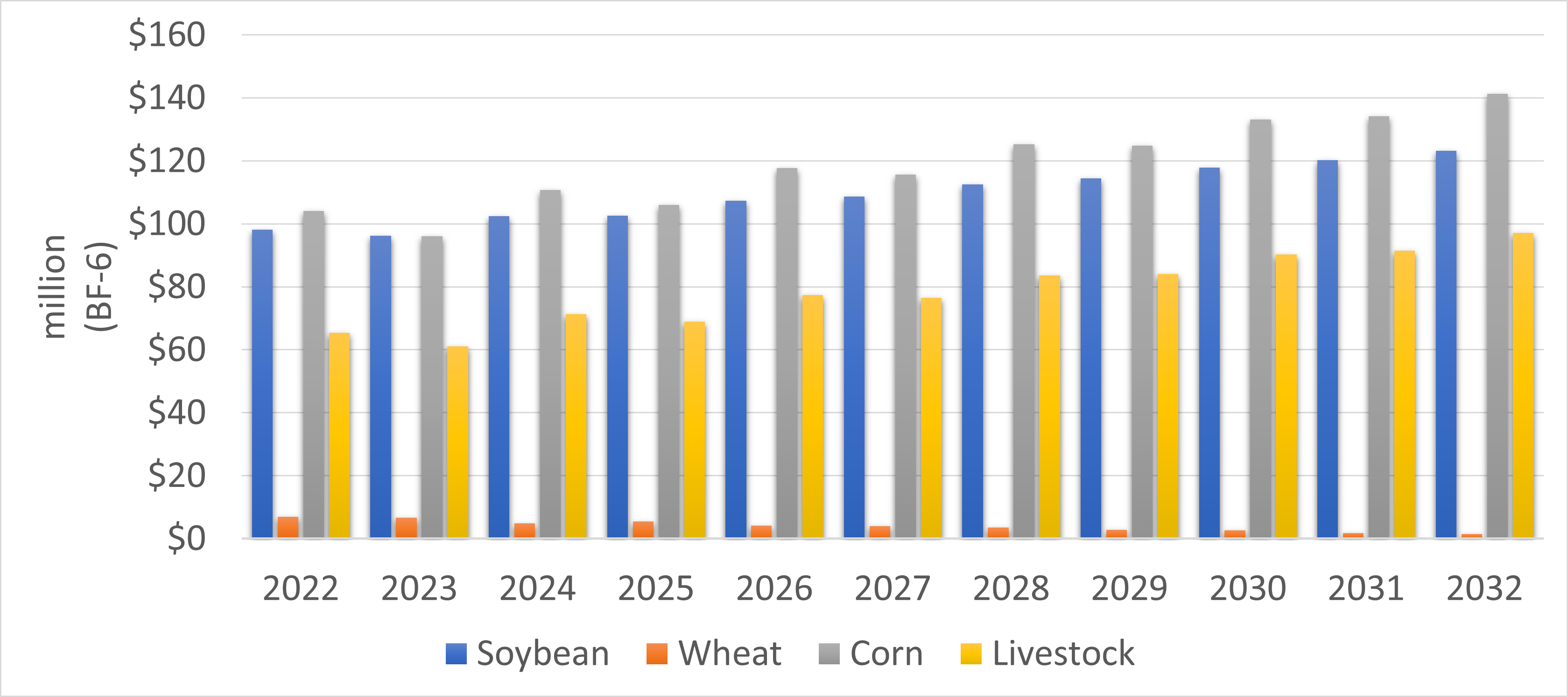}
    \caption{Cost savings in agriculture (BF-6; see Section \ref{sec 4.6}).}
    \label{fig6}
\end{figure}

Investment in AAM infrastructure can be justified in terms of the significant benefits gained by deploying AAM. The yearly estimated values of travel time savings for passengers (BF-1), passengers' yearly safety cost reduction (BF-2) and cost and travel time delay savings for bridge inspections (BF-5) are depicted in Figure \ref{fig4}. The figure shows an increasing trend in travel time cost saving over the analysis period. Ohio passengers' travel time cost saving in 2032 is estimated to be over \$22 million, with a 65.2\% CAGR during 2022 to 2032 (see Equation (\ref{12}) for the calculation of CAGR). It is projected to be valued over \$0.7 million by 2032, growing at a 70.5\% CAGR for the next 11 years. Based on these estimates, AAM infrastructure has the potential to provide a safer transportation system for passengers relative to the existing car-based ground transportation system. An estimated cost reduction of \$8.78 million can be achieved by the government from the starting year of our analysis by utilizing drones for inspections. It was estimated to generate over \$9.75 million worth of cost reduction in bridge inspections by 2032, following an increasing linear trend as shown in Figure \ref{fig4}. These estimates are also an indication of the opportunities present in other infrastructure inspections---such as tunnel, pipeline, highway, and tower inspections---which are more or less similar to bridge inspections.

The yearly estimated benefit values associated with cost and delivery time savings in package delivery by drones (BF-3), and time and inventory cost savings in cargo delivery by eVTOLs (BF-4) are illustrated in Figure \ref{fig5}. The projected exponential trend of BF-3 suggests that AAM will fuel the growth of package delivery time savings. By 2032, the value of cost and delivery time savings is estimated to reach \$7.8 billion, exhibiting a CAGR of 60.2\% for the analysis period, in response to the growing demand for quick delivery. The results indicate that Ohio's cargo delivery by eVTOLs will provide substantial lead time savings and cost reductions for businesses that need to deliver products to consumers over the next 11 years, giving another strong reason to invest in AAM infrastructure. More than \$66 billion valuation of BF-4 by 2032, at a 63.36\% CAGR during the forecasted period, indicates that the AAM infrastructure can generate meaningful economic activity in the logistic network by providing the scopes for more effective use of warehouse space along with increasing the lead time savings.

\begin{figure}[H]
    \includegraphics[width=13cm,height=5cm]{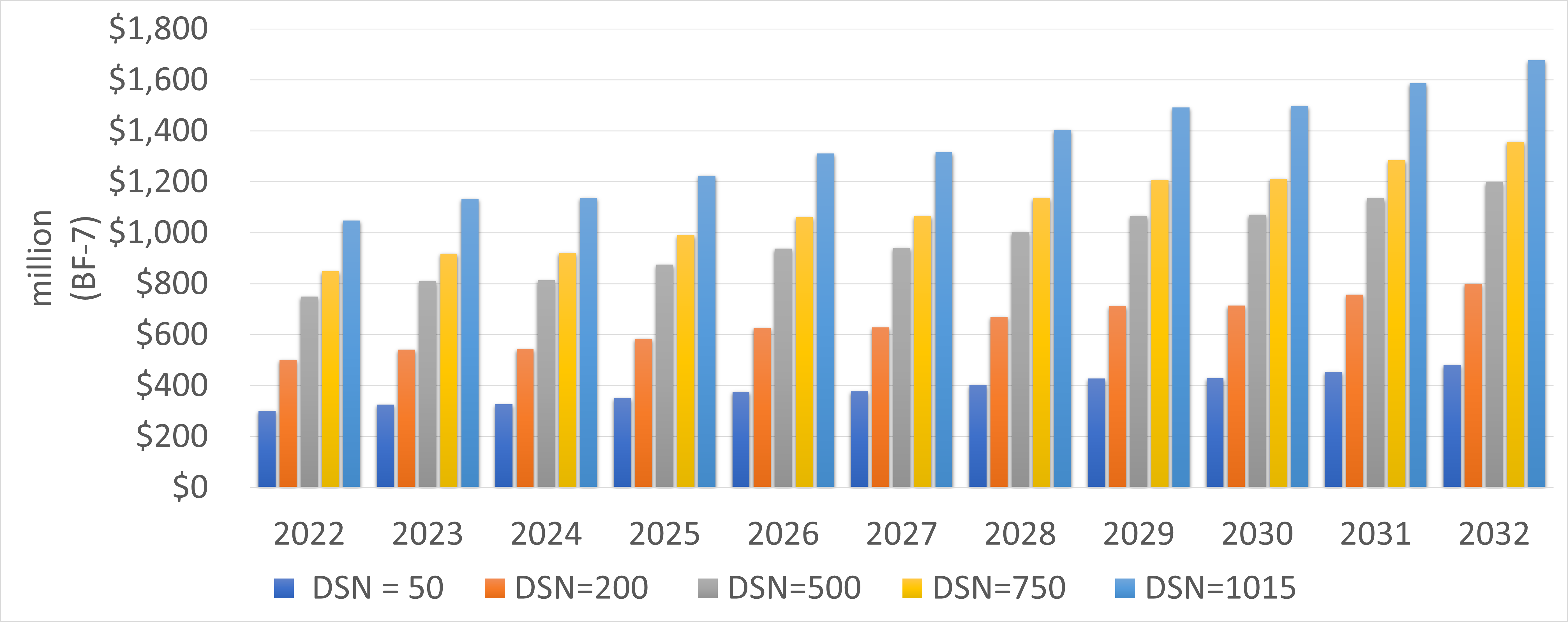}
    \caption{Life-Saving value in AED Delivery (BF-7; see Section \ref{sec 4.7}).}
    \label{fig7}
\end{figure}

\vspace{-6pt}
\begin{figure}[H]
    \includegraphics[width=13cm,height=6cm]{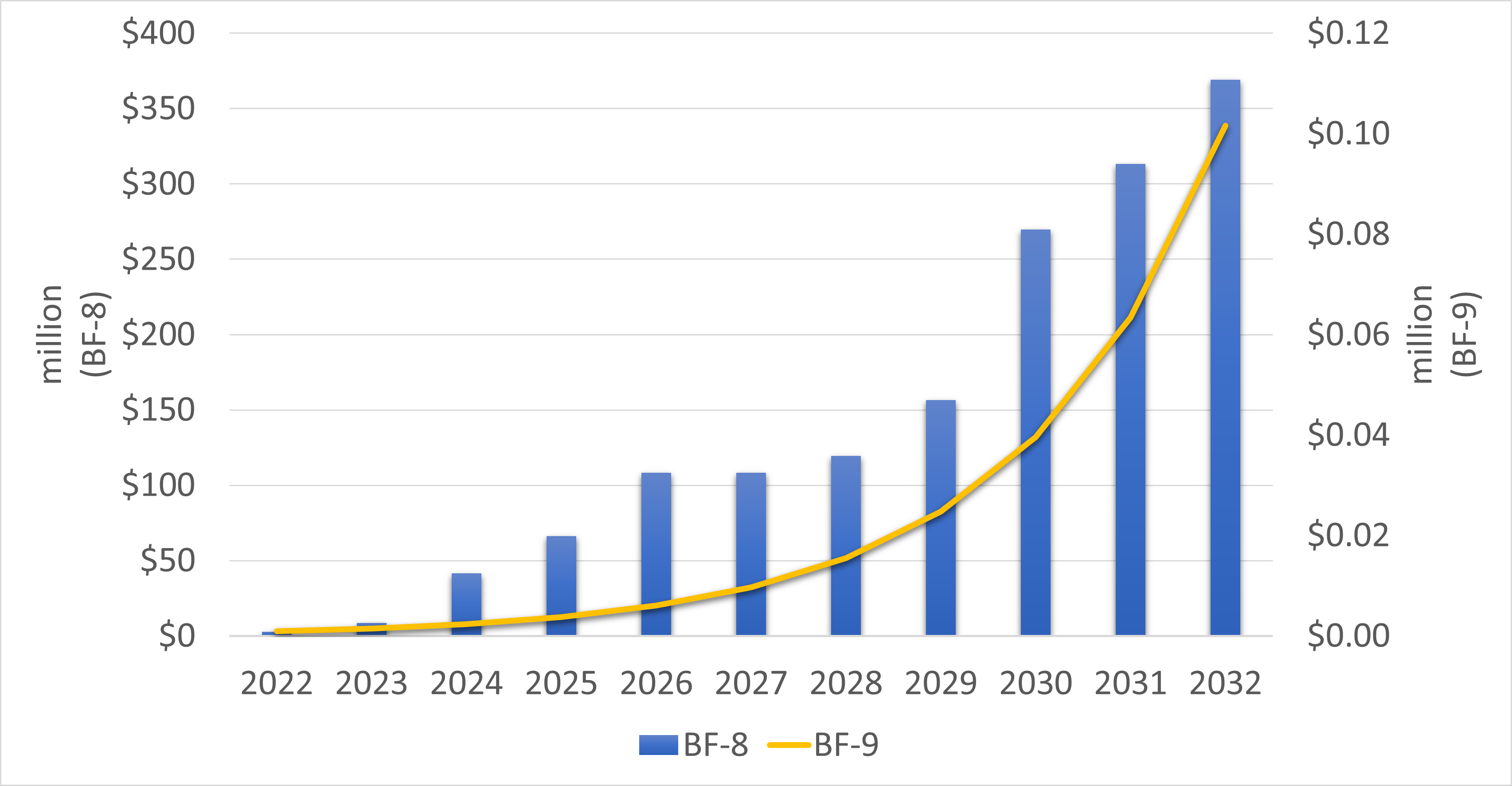}
    \caption{Yearly estimated values coming from BF-8 (see Section \ref{sec 4.8}) and BF-9 (see Section \ref{sec 4.9}).}
    \label{fig8}
\end{figure}

\vspace{-6pt}
\begin{figure}[H]
    \includegraphics[width=13cm,height=6cm]{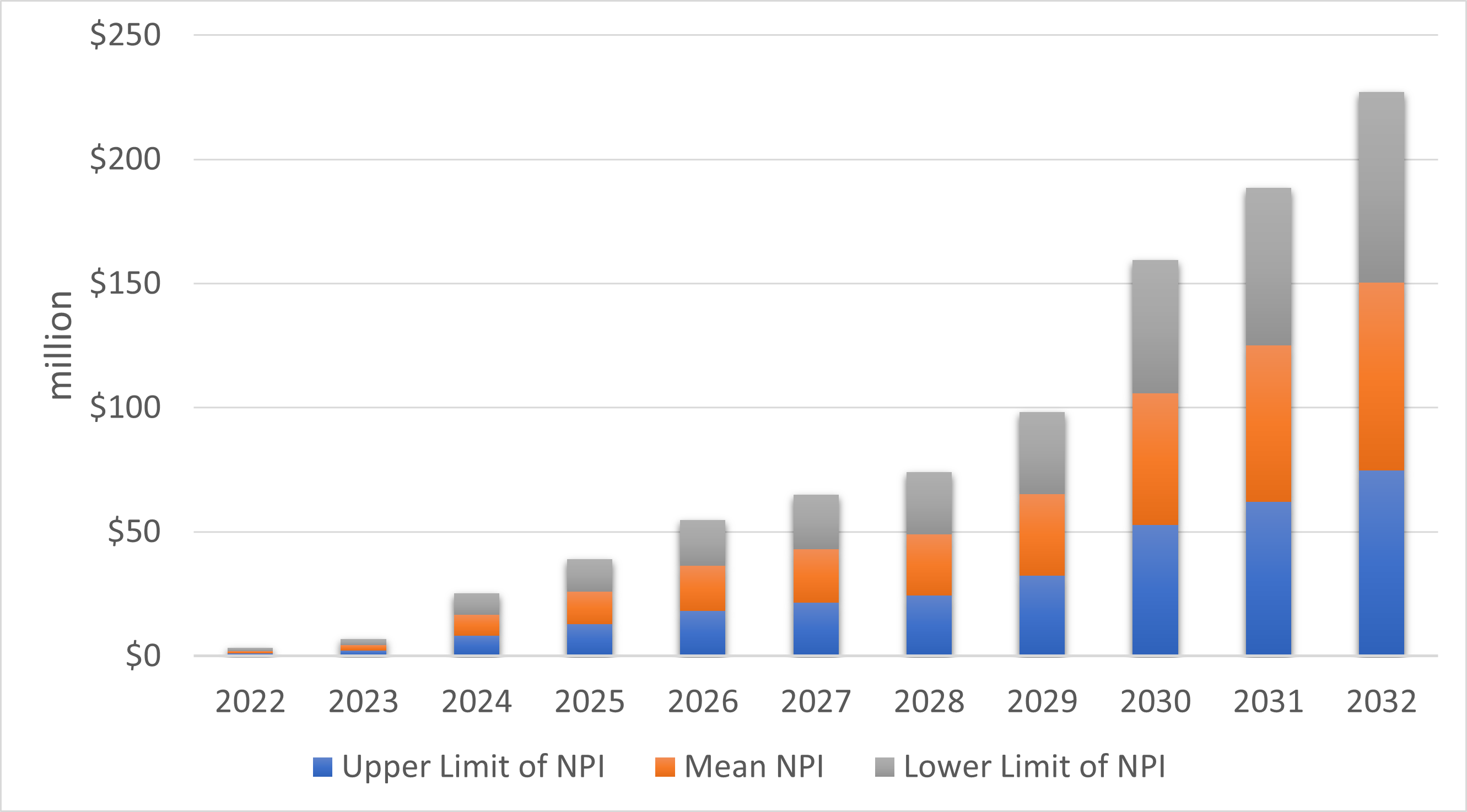}
    \caption{Estimated yearly net positive gain considering nine benefit factors and CAPEX-OPEX cost.}
    \label{fig9}
\end{figure}

Further justification for building AAM infrastructure in Ohio can be seen in Figure \ref{fig6}, which shows an estimated \$274 million of cost reduction to be experienced by the farmers in the crop and livestock farming operations at the beginning of our analysis span. By 2032, this number will exceed \$363 million. Focusing on the values for the individual crops, it can be seen that the savings will be higher for corn and soybean relative to wheat. Though we studied only the OHCA case in this paper, Figure \ref{fig7} helps to realize the potential of AAM infrastructure in improving the medical delivery system. For the year 2032, the projected life-saving value derived from using drones for AED delivery for patients is around \$481 million when DSN (see Section \ref{sec 4.7}) is 50, which increases to \$1677 million when DSN is 1015. During the evaluation of this benefit factor, we adjusted the life-saving value to take into account the cost needed for renovating the local police, fire, or emergency stations to install the necessary drone docking stations. As shown in the figure, the life-saving value is going higher while DSN is increasing. Thus, it can be interpreted that the higher necessity of response to the medical emergency cases can offset the capital and operating cost of DSN.

Next, Figure \ref{fig8} indicates increased tax income for government (BF-8) and savings in social cost of greenhouse gases (BF-9). An upward trend of estimated tax income is projected for the state government of Ohio coming from services utilizing the AAM infrastructure. Over the 11-year time period, it is estimated to reach \$369 million by 2032 at a 63.3\% CAGR. An exponential upward trend was also observed for BF-9, with around \$0.1 million worth of reduction of social costs of GHG in 2032 at 59.8\% CAGR. These estimates suggest that AAM infrastructure has significant potential in contributing towards sustainability in transportation relative to the traditional ground transportation. 

Finally, Figure \ref{fig9} displays the net positive gain for society and environment for the next 11 years generated from AAM infrastructure in Ohio. For 2022, the estimated net positive impact (NPI) is \$9.74 billion, and following 61.6\% CAGR, the gain is estimated to rise at \$76.68~billion. We calculated this gain by subtracting the CAPEX-OPEX cost (see Figure \ref{fig3}) from the total sum of all the values of the nine benefit factors over the 11 years. For addressing uncertainty in forecasted variables, we considered 95\% confidence interval during forecasting each variable and presented the outcome by upper limit, mean and lower limit of NPI in Figure \ref{fig9}.

\section{Conclusions}\label{sec7}
In this paper, we performed a cost--benefit analysis of AAM infrastructure from the perspective of the government, primarily focusing on AAM's non-monetary benefits for the society and environment. The state of Ohio was taken as our case study. Nine different types of primary benefits of AAM were identified, which mostly arise from the various use cases enabled and supported by AAM. These benefits were quantified in dollar amounts for the 2022--2032 time period. Significant cost reductions and time savings can be achieved in these use cases by AAM compared to surface transportation and other traditional approaches. The benefits were also found to substantially surpass the costs associated with AAM infrastructure. Our findings provide strong justification for investment by government in AAM infrastructure in Ohio, and, by extension, in other states and regions with similar characteristics.

In this competitive world, time savings and cost reductions are important terms to make a strategic decision before investing on a project. The significant amount of savings projected in this paper will help the decision-makers to take initiative towards AAM implementation. Among the nine benefit factors, time and inventory cost savings in cargo delivery by eVTOLs (BF-4), and cost and delivery time savings in package delivery by drones (BF-3) have the largest gain with \$66.76 billion and \$7.82 billion, respectively. These values suggest that the AAM implementation has the potential to produce a considerable economic impact for the logistic companies, which will allow more efficiency for handling warehouse space and lead time. Other benefit factors also present notable insights to make decision---such as increased crop yield can address national and global food shortages and lower food prices, decreased carbon footprint can reduce the issue related to global warming, reduced delivery time in response to medical emergency case like OHCA can increase the survival rate, etc. Understanding these advantages can aid in gaining community acceptance of AAM, which is crucial for its effective implementation.

There are several interesting expansions to this work that make for compelling studies in the future. Firstly, the AAM cost--benefit identification, quantification and analysis approach used in this paper can be applied to other potential states of the USA to assess the impact of AAM in those places before any investment is made there. In our study, we considered one medical emergency case. Secondly, in future studies, the scope of this work can be expanded to include benefits derived other envisioned use cases of AAM---such as organ, isotope, vaccine, blood, lab specimen delivery, search and rescue people in potential danger, or vessels lost at sea, security (police) patrol and law enforcement, mapping and surveying of (inaccessible) terrains and locations, aerial streaming of live events, firefighting, meteorology measurements, etc. In addition, tunnel, pipeline, tower, highway, and other infrastructure inspections can be included in future work. Lastly, future work can investigate effects of AAM which can have potentially negative impacts on the society---such as, increase in air traffic congestion, flight delays and noise levels; impact on ground vehicle drivers due to visual distraction by AAM aircraft; and privacy invasion.

\vspace{6pt}

\funding{This work was supported by the Ohio Department of Transportation (Agreement \mbox{No.: 36496}, PID: 114242, SJN: 136337). Any opinions, findings, and conclusions or recommendations expressed in this material are those of the authors and do not necessarily reflect the views of the Ohio Department of Transportation. The authors would like to thank Rubén Del Rosario and Shahab Hasan at Crown Consulting for sharing relevant data from their prior AAM research.}

\dataavailability{The data presented in this study are available on request from the corresponding author. The data are not publicly available due to privacy restriction.
} 

\conflictsofinterest{The authors declare no conflict of interest.}

\appendixtitles{no} 
\appendixstart
\appendix
\section[\appendixname~\thesection]{\label{appendix}}
Tables \ref{table3} and \ref{tab4} are showing the operational cost breakdown for a bridge inspection using snooper and drone, respectively. The data were taken from \cite{flow1}.

\begin{table}
\caption{Operational cost rate per bridge inspection using snooper (source: \cite{flow1}). \label{table3}}
\begin{tabular}{m{1cm}m{1cm}m{4.2cm}}
\toprule
\multicolumn{3}{c}{\textbf{Labor}} \\ \midrule
\multicolumn{1}{c}{Employee} & \multicolumn{1}{c}{Bridge specialist} & Highway technician \\ \midrule
\multicolumn{1}{c}{No. of   employee} & \multicolumn{1}{c}{3} & \multicolumn{1}{c}{3} \\ \midrule
\multicolumn{1}{c}{Hours per   employee} & \multicolumn{1}{c}{8} & \multicolumn{1}{c}{8} \\ \midrule
\multicolumn{1}{c}{Hourly rate} & \multicolumn{1}{c}{\$37} & \multicolumn{1}{c}{\$21} \\ \midrule
\multicolumn{1}{c}{Fringe} & \multicolumn{1}{c}{45\%} & \multicolumn{1}{c}{45\%} \\ \midrule
\multicolumn{1}{c}{Total cost} & \multicolumn{1}{c}{\$854} & \multicolumn{1}{c}{\$727} \\ \midrule
\multicolumn{1}{c}{Total   payroll per bridge inspection, A} & \multicolumn{2}{c}{\$2018} \\ \midrule
\multicolumn{3}{c}{\textbf{Equipment}} \\ \midrule
\multicolumn{1}{c}{Equipment type} & \multicolumn{1}{c}{Snooper truck} &Other equipment   (pickups, cones, sign trucks, etc.) \\ \midrule
\multicolumn{1}{c}{Equipment   rate} & \multicolumn{1}{c}{\$625} & \multicolumn{1}{c}{\$500} \\ \midrule
\multicolumn{1}{c}{Total equipment cost, B} & \multicolumn{2}{c}{\$1125} \\ \midrule
\multicolumn{3}{c}{\textbf{Cost rate per bridge inspection}} \\ \midrule
\multicolumn{1}{c}{Cost rate (A+B) based on core hours} & \multicolumn{2}{c}{\$3143} \\ \midrule
\multicolumn{1}{c}{Cost rate based on nights/weekends} & \multicolumn{2}{c}{\$4152} \\ \bottomrule                 
\end{tabular}
\end{table}

\begin{table}[H]
\caption{Operational cost rate per bridge inspection using drone (source: \cite{flow1}). \label{tab4}}
\tabcolsep=0.7cm
\begin{tabular}{m{6cm}m{10cm}}
\toprule
\multicolumn{2}{c}{\textbf{Labor}} \\ \midrule
\multicolumn{1}{c}{Employee} & \multicolumn{1}{c}{Bridge specialist} \\ \midrule
\multicolumn{1}{c}{No. of  employee} & \multicolumn{1}{c}{2} \\ \midrule
\multicolumn{1}{c}{Hours per employee} & \multicolumn{1}{c}{4} \\ \midrule
\multicolumn{1}{c}{Hourly rate} & \multicolumn{1}{c}{\$37} \\ \midrule
\multicolumn{1}{c}{Fringe} & \multicolumn{1}{c}{45\%} \\ \midrule
\multicolumn{1}{c}{Total payroll per bridge inspection, A} & \multicolumn{1}{c}{\$427} \\ \midrule
\multicolumn{2}{c}{\textbf{Equipment}} \\ \midrule
\multicolumn{1}{c}{Equipment type} & \multicolumn{1}{c}{Drone, software, tablet \& pickups} \\ \midrule
\multicolumn{1}{c}{Total equipment cost rate, B} & \multicolumn{1}{c}{\$95} \\ \midrule
\multicolumn{2}{c}{\textbf{Cost rate per bridge inspection}} \\ \midrule
\multicolumn{1}{c}{Cost rate (A+B) based on core hours} & \multicolumn{1}{c}{\$522} \\ \midrule
\multicolumn{1}{c}{Cost rate based on nights/weekends} & \multicolumn{1}{c}{\$735} \\ \bottomrule
\end{tabular}
\end{table}

\begin{adjustwidth}{-\extralength}{0cm}

\reftitle{References}

\end{adjustwidth}

\end{document}